\documentclass[12pt]{article}
\usepackage{amsmath,amssymb}
\usepackage{graphicx}
\newcommand{\eqdef}{\stackrel{\text{def}}{=}}

\newcommand{\n}{\nonumber\\}

\newcommand{\ignore}[1]{}
\numberwithin{equation}{section}
\newcommand{\Romannumeral}[1]{\uppercase\expandafter{\romannumeral#1}}

\newtheorem{theo}{\bf Theorem}[section]
\newtheorem{lemma}[theo]{\bf Lemma}
\newtheorem{coro}[theo]{\bf Corollary}
\newtheorem{rema}[theo]{\bf Remark}

\allowdisplaybreaks[4]

\begin{document}

\baselineskip=20pt

%%%%%%%%%%%%%%%%%%%%%%%%%%%%%%%%%%%%%%%%%%%%%%%%%%%%%%%%%%%%
%                                                          %
%  Title page                                              %
%                                                          %
%%%%%%%%%%%%%%%%%%%%%%%%%%%%%%%%%%%%%%%%%%%%%%%%%%%%%%%%%%%%
\newcommand{\preprint}{
\vspace*{-20mm}%\begin{flushleft} Version:  February 18 2014 \end{flushleft}
%    \begin{flushright}\normalsize \sf
%     DPSU-13-3\\
%     {\tt arXiv:1306.nnnn[math-ph]}\\
%     July 2013
%   \end{flushright}
}
\newcommand{\Title}[1]{{\baselineskip=26pt
   \begin{center} \Large \bf #1 \\ \ \\ \end{center}}}
\newcommand{\Author}{\begin{center}
   \large \bf  Ryu Sasaki \end{center}}
\newcommand{\Address}{\begin{center}
      Department of Physics, Shinshu University,\\
     Matsumoto 390-8621, Japan\\
     %and\\
     Department of Physics, National Taiwan University, \\
     Taipei 10617, Taiwan
    \end{center}}
\newcommand{\Accepted}[1]{\begin{center}
   {\large \sf #1}\\ \vspace{1mm}{\small \sf Accepted for Publication}
   \end{center}}

\preprint
\thispagestyle{empty}

\Title{Perturbations around the zeros of classical orthogonal polynomials}                     %

\Author

\Address
%\vspace{1cm}

\begin{abstract}
Starting from degree $\mathcal{N}$ solutions of a time dependent Schr\"odinger-like
equation for classical orthogonal polynomials, a linear matrix equation  describing 
perturbations around the $\mathcal{N}$ zeros of the polynomial is derived.
The matrix has remarkable Diophantine properties. Its eigenvalues are independent
of the zeros.
The corresponding eigenvectors provide the representations of the lower degree 
($0,1,\ldots,\mathcal{N}-1$) polynomials in terms of the zeros 
of the degree $\mathcal{N}$
polynomial. The results are valid universally for all the classical orthogonal polynomials,
including the Askey scheme of hypergeometric orthogonal polynomials and its $q$-analogues.
\end{abstract}

%%%%%%%%%%%%%%%%%%%%%%%%%%%%%%%%%%%%%%%%%%%%%%%%%%%%%%%%%%%%%%%
%                                                             %
%  1. Introduction                                            %
%                                                             %
%%%%%%%%%%%%%%%%%%%%%%%%%%%%%%%%%%%%%%%%%%%%%%%%%%%%%%%%%%%%%%%
\section{Introduction}
\label{intro}

The properties of the zeros of orthogonal polynomials, in particular, of
classical orthogonal polynomials have been a fascinating subject for many years 
\cite{stiel, szego}. In this paper we define {\em classical orthogonal polynomials\/}
 as  polynomials {\em 
satisfying second order differential or difference equations and the three term recurrence relations\/} \cite{chihara}.
That is, the Askey scheme of hypergeometric orthogonal polynomials 
and its $q$-analogues are included \cite{askey,ismail,koekswart} 
but not the recently discovered 
 multi-indexed \cite{os25} and  exceptional \cite{xop} orthogonal polynomials.
Here we report a small contribution to the subject by presenting an 
$\mathcal{N}\times \mathcal{N}$ matrix $\mathcal{M}$ \eqref{Mdef} possessing 
remarkable Diophantine properties.
The matrix $\mathcal{M}$ describes the small oscillations around the zeros of a degree
$\mathcal{N}$ classical orthogonal polynomial.
Main results are (i) The eigenvalues of $\mathcal{M}$ are independent of the zeros.
They are the difference of those of the
differential/difference operator $\widetilde{\mathcal H}$ \eqref{Htil},
which governs the classical orthogonal polynomial, corresponding to the degree 
$\mathcal{N}$ and a lower degree, Theorem\,\ref{theo:1}.
(ii) The corresponding eigenvectors provide {\em representations of the lower degree
polynomials in terms of the zeros of
degree ${\mathcal N}$ polynomial\/}, Theorem\,\ref{theo:2}.
The theorems are universal, meaning that they apply to all the classical orthogonal polynomials.

The idea of the present research was influenced by a recent paper of Bihun and Calogero
\cite{bihun-cal} discussing the multi-particle interactions and their equilibrium
implied by the classical orthogonal polynomials, in particular, the Wilson and the Racah polynomials.

This paper is organized as follows.
In section two, we start with the time dependent Schr\"odinger equation \eqref{scheq} 
of quantum 
mechanics, which offers a unified framework for discussing the classical orthogonal
polynomials as the main parts of the eigenfunctions.
By restricting the general solution \eqref{gensol} to those up to degree $\mathcal{N}$,
a time dependent equation \eqref {tevo} for degree $\mathcal{N}$ 
polynomials is obtained
in \S\ref{sec:poly}. Section \ref{sec:pert} is the main part of the paper.
By considering perturbations around the zeros of the degree $\mathcal{N}$
classical polynomials, the time-dependent equation is rewritten as a linear equation  
with an
$\mathcal{N}\times\mathcal{N}$ matrix $\mathcal{M}$ describing the infinitesimal oscillations around the zeros.
The main Theorems about the eigenvalues and eigenfunctions are stated as 
the natural consequence of the construction.
In sections four, five and six, the various data 
of the classical orthogonal polynomials are
presented for the explicit construction of the matrix $\mathcal{M}$ for verification.
Section four is for the Classical orthogonal polynomials, that is, the Hermite,
Laguerre and Jacobi.
Section five is for the ($q$-)Askey scheme of hypergeometric polynomials having 
the continuous orthogonality weight functions, that is the Wilson 
and Askey-Wilson polynomials and their reduced form polynomials.
Section six is for the ($q$-)Askey scheme of hypergeometric polynomials having 
the purely discrete weight functions, that is the Racah and $q$-Racah polynomials 
and their reduced form polynomials. The final section is for a summary and comments.
A small Appendix is for the symbols and definitions
related to the ($q$-)hypergeometric functions.

%%%%%%%%%%%%%%%%%%%%%%%%%%%%%%%%%%%%%%%%%%%%%%%%%%%%%%%%%%%%%%%
%                                                             %
%  2. Time dependent Schr\"odinger equations              %
%                                                             %
%%%%%%%%%%%%%%%%%%%%%%%%%%%%%%%%%%%%%%%%%%%%%%%%%%%%%%%%%%%%%%
\section{Time dependent Schr\"odinger equations}
\label{sec:eigen}

The time dependent Schr\"odinger equation
\begin{equation}
i\frac{\partial \Psi(x,t)}{\partial t}=\mathcal{H} \Psi(x,t)
\label{scheq}
\end{equation}
is exactly solvable if the corresponding eigenvalue problem of the Hamiltonian or the
Schr\"odinger operator
$\mathcal{H}$
\begin{equation}
\mathcal{H}\phi_n(x)=\mathcal{E}(n)\phi_n(x),\quad n=0,1,\ldots,
\label{eigpro}
\end{equation}
is exactly solvable.
Throughout this paper the Hamiltonian $\mathcal{H}$ is assumed to be time independent.
In terms of the complete set of solutions $\{\mathcal{E}(n),\phi_n(x)\}$ of 
the eigenvalue problem \eqref{eigpro}, the general solution of the time dependent
Schr\"odinger equation is given by
\begin{equation}
\Psi(x,t)=\sum_{n=0}^\infty c_n e^{-i\mathcal{E}(n)t}\phi_n(x),
\label{gensol}
\end{equation}
in which $\{c_n\}$ are the constants of integration.

Hereafter we discuss one dimensional exactly solvable Hamiltonian systems 
\cite{infhull,susyqm,os12,os13,os24}, with a Hamiltonian $\mathcal{H}$ which is
a second order differential or difference operator.
Their eigenfunctions
include {\em classical orthogonal polynomials\/}, {\em i.e.},
the Hermite, Laguerre, Jacobi, Wilson, Askey-Wilson, Racah and $q$-Racah polynomials, 
and others \cite{os25, xop}.
In other words, all the hypergeometric orthogonal polynomials of Askey scheme constitute
the main part of the eigenfunctions of certain `discrete' 
quantum mechanical systems with
pure imaginary and/or real shifts \cite{os12,os13,os24}.
For most solvable examples, the eigenfunctions have a factorised form
\begin{equation}
\phi_n(x)=\phi_0(x)P_n\left(\eta(x)\right).
\label{fac}
\end{equation}
Here $\phi_0(x)$ is the ground state  wave function 
and its square $\phi_0(x)^2$ provides
the orthogonality weight function for the polynomial $P_n\left(\eta(x)\right)$
 of degree $n$%
 \footnote{In the case of recently discovered {\em multi-indexed\/} 
 orthogonal polynomials \cite{os25}, 
 $n$ stands for the number of nodes in the orthogonality interval.
 The degree of the polynomial is greater than $n$.
 }
  in $\eta(x)$,
which is called the {\em sinusoidal coordinate} \cite{os7}.

By similarity transforming
the Hamiltonian in terms of the ground state wave function $\phi_0(x)$, we obtain 
the differential/difference operator $\widetilde{\mathcal H}$:
\begin{equation}
\widetilde{\mathcal H}\eqdef \phi_0(x)^{-1}\circ\mathcal{H}\circ\phi_0(x),
\label{Htil}
\end{equation}
governing the classical polynomials $\{P_n\left(\eta(x)\right)\}$:
\begin{equation}
\widetilde{\mathcal H}P_n\left(\eta(x)\right)
=\mathcal{E}(n)P_n\left(\eta(x)\right),\quad
n=0,1,\ldots .
\label{HtilPn}
\end{equation}
In other words,  
$\widetilde{\mathcal H}$ keeps the polynomial space 
$\{1,P_1(\eta(x)),P_2(\eta(x)),\ldots,P_n(\eta(x))\}$ invariant.

In the rest of this paper, we consider only those systems having the above factorised
eigenfunctions \eqref{fac}.
We further restrict our attention to the {\em classical orthogonal polynomials\/} 
$P_n\bigl(\eta(x)\bigr)$, that is, the new (the multi-indexed \cite{os25} and
exceptional \cite{xop})
orthogonal polynomials will not be included.

%%%%%%%%%%%%%%%%%%%%%%%%%%%%%%%%%%%%%%%%%%%%%
%                                           %
% 2.1 Polynomial solutions                    %
%                                           %
%%%%%%%%%%%%%%%%%%%%%%%%%%%%%%%%%%%%%%%%%%%%%
\subsection{Polynomial solutions}
\label{sec:poly}

Let us fix a positive integer $\mathcal{N}$ 
and restrict the general solution \eqref{gensol}
to those having {\em degrees up to \/} $\mathcal{N}$:
\begin{align} 
 \Psi_{\mathcal N}(x,t)&=\sum_{n=0}^{\mathcal N} 
 c_n e^{-i\mathcal{E}(n)t}\phi_n(x),\\
 &=e^{-i\mathcal{E}(\mathcal{N})t}\phi_0(x)\psi_{\mathcal N}(x,t).
 \label{polysol}
\end{align}
Here the function
\begin{equation}
\psi_{\mathcal N}(x,t)\eqdef \sum_{n=0}^{\mathcal N}c_n 
e^{i(\mathcal{E}(\mathcal{N})-\mathcal{E}(n))t}P_n\left(\eta(x)\right)
\label{psiNdef}
\end{equation}
is a polynomial of degree $\mathcal{N}$ in $\eta(x)$.
We choose the coefficient $c_{\mathcal N}$ of the highest degree polynomial  
$P_{\mathcal N}\left(\eta(x)\right)$
to make it monic:
\begin{equation}
c_{\mathcal N}P_{\mathcal N}\bigl(\eta(x)\bigr)
=\prod_{n=1}^{\mathcal N}\left(\eta(x)-\eta(x_n)\right),
\label{xndef}
\end{equation}
in which $\{\eta(x_n)\}$, $n=1,\ldots,\mathcal{N}$,  
are the {\em zeros} of $P_{\mathcal N}\left(\eta(x)\right)$.
The polynomial $\psi_{\mathcal N}(x,t)$ satisfies the time evolution equation
\begin{align}
\frac{\partial \psi_{\mathcal N}(x,t)}{\partial t}
=-i\widetilde{\mathcal H}_{\mathcal N}\psi_{\mathcal N}(x,t),\qquad
\widetilde{\mathcal H}_{\mathcal N}\eqdef 
\widetilde{\mathcal H}-\mathcal{E}({\mathcal N}).
\label{tevo}
\end{align}
For a given set of parameters $\{c_n\}$, $n=0,\ldots,\mathcal{N}-1$, the polynomial 
$\psi_{\mathcal N}(x,t)$ can be regarded as a $t$-dependent deformation
of the highest degree monic polynomial $c_{\mathcal N}P_{\mathcal N}\left(\eta(x)\right)$:
\begin{equation}
\psi_{\mathcal N}(x,t)=\prod_{n=1}^{\mathcal N}\bigl(\eta(x)-\eta(x_n(t))\bigr),
\label{xntdef}
\end{equation}
in which $\{x_n(t)\}$ are certain $t$-dependent functions, describing the zeros of 
$\psi_{\mathcal N}(x,t)$ at time $t$.

%%%%%%%%%%%%%%%%%%%%%%%%%%%%%%%%%%%%%%%%%%%%%
%                                           %
% 3 Perturbations around the zeros                   %
%                                           %
%%%%%%%%%%%%%%%%%%%%%%%%%%%%%%%%%%%%%%%%%%%%%
\section{Perturbations around the zeros}
\label{sec:pert}

Among the generic $t$-dependent deformations \eqref{xntdef} of the classical polynomial
$P_{\mathcal N}\left(\eta(x)\right)$, let us focus on those describing 
{\em infinitesimal
oscillations around the zeros\/} of $P_{\mathcal N}\left(\eta(x)\right)$:
\begin{equation}
x_n(t)=x_n+\epsilon\gamma_n(t),\quad 0<\epsilon\ll1,\quad n=1,\ldots,\mathcal{N}.
\label{ansatz}
\end{equation}
In other words, instead of the general deformation
\eqref{xntdef} by $\{c_n\}$, we choose infinitesimal $\{c_n\}$ 
so that the deformation can be considered as perturbations around the 
zeros of $P_{\mathcal N}(\eta(x))$.
The above ansatz \eqref{ansatz} leads to
\begin{align} 
\psi_{\mathcal N}(x,t)&=\prod_{n=1}^{\mathcal N}\bigl(\eta(x)-\eta(x_n)\bigr)
             -\epsilon\sum_{n=1}^{\mathcal N}\gamma_n(t)\dot{\eta}(x_n)
    \prod_{j\neq n}^{\mathcal N}\bigl(\eta(x)-\eta(x_j)\bigr) +O(\epsilon^2),\\
 \text{with} \quad \dot{\eta}(x)&\eqdef \frac{d\eta(x)}{dx}.
\end{align}
With this ansatz, the l.h.s. of the time evolution equation \eqref{tevo} 
is a degree $\mathcal{N}-1$ polynomial in $\eta(x)$:
\begin{equation}
-\epsilon\sum_{n=1}^{\mathcal N}\frac{d\gamma_n(t)}{dt}\dot{\eta}(x_n)
             \prod_{j\neq n}^{\mathcal N}\bigl(\eta(x)-\eta(x_j)\bigr) +O(\epsilon^2).
\end{equation}
The r.h.s. is also a degree $\mathcal{N}-1$ polynomial in $\eta(x)$:
\begin{equation}
i\epsilon\sum_{m=1}^{\mathcal N}\gamma_m(t)\dot{\eta}(x_m)
             \widetilde{\mathcal H}_{\mathcal N}\prod_{j\neq m}^{\mathcal N}
             \bigl(\eta(x)-\eta(x_j)\bigr) +O(\epsilon^2),
\end{equation}
since the leading polynomial $P_{\mathcal N}\left(\eta(x)\right)\propto 
\prod_{n=1}^{\mathcal N}\left(\eta(x)-\eta(x_n)\right)$ is annihilated by 
$\widetilde{\mathcal H}_{\mathcal N}$:
\begin{equation}
\widetilde{\mathcal H}_{\mathcal N}P_{\mathcal N}\left(\eta(x)\right)=0.
\label{HtilN}
\end{equation}

The polynomial evolution equation \eqref{tevo}, being of an $\mathcal{N}-1$  degree,
 is satisfied when its evaluation at $\mathcal{N}$ independent points are satisfied
 and vice versa.
 Without loss of generality, we can choose the $\mathcal{N}$ zeros  $\{\eta(x_n)\}$ of 
 $P_{\mathcal N}\left(\eta(x)\right)$.
This leads to $\mathcal{N}$ linear ODE's for the unknown functions $\{\gamma_n(t)\}$
at the leading order of $\epsilon$:
\begin{align} 
&\frac{d\gamma_n(t)}{dt}\dot{\eta}(x_n) 
\prod_{j\neq n}^{\mathcal N}\bigl(\eta(x_n)-\eta(x_j)\bigr)\n
\quad &=-i\sum_{m=1}^{\mathcal N}\gamma_m(t)\dot{\eta}(x_m)
\left.\left(\widetilde{\mathcal H}_{\mathcal N}
\prod_{j\neq m}^{\mathcal N}\bigl(\eta(x)-\eta(x_j)\bigr)
\right)\right |_{x=x_n},\quad n=1,\ldots,\mathcal{N},
\end{align}
which can be rewritten in a matrix form:
\begin{align}   
\frac{d\gamma_n(t)}{dt}&= i\sum_{m=1}^{\mathcal N}\mathcal{M}_{n\,m}\gamma_m(t), 
\quad  n=1,\ldots,\mathcal{N}, 
\label{mateq}\\
\mathcal{M}_{n\,m}&\eqdef -\frac{\dot{\eta}(x_m)\left.
\left(\widetilde{\mathcal H}_{\mathcal N}
\prod_{j\neq m}^{\mathcal N}\bigl(\eta(x)-\eta(x_j)\bigr)
\right)\right |_{x=x_n}}{\dot{\eta}(x_n)
\prod_{j\neq n}^{\mathcal N}\bigl(\eta(x_n)-\eta(x_j)\bigr)}.
\label{Mdef}
\end{align}

By construction, we have the following
\begin{theo}\label{theo:1}
The eigenvalues of $\mathcal{M}$ are
\begin{equation}
\mathcal{E}(\mathcal{N})-\mathcal{E}(m),\quad m=0,1,\ldots,\mathcal{N}-1,
\label{eigvals}
\end{equation}
which depend on the basic parameters of $\widetilde{\mathcal H}$ but do not depend on the zeros $\{\eta(x_n)\}$ directly.
\end{theo}
\begin{theo}\label{theo:2}
The corresponding eigenvectors $\{{\rm v}^{(m)}_n\}$ of $\mathcal{M}$, 
\begin{equation}
\sum_{\ell=1}^{\mathcal N}\mathcal{M}_{n\,\ell}{\rm v}^{(m)}_{\ell}
=(\mathcal{E}(\mathcal{N})-\mathcal{E}(m))
{\rm v}^{(m)}_n
\label{eigveceq}
\end{equation}
yield the representations of the lower degree polynomials $\{P_m(\eta)\}$, 
$m=0,1,\ldots,\mathcal{N}-1$, in terms of the zeros 
$\{\eta(x_n)\}$ of $P_{\mathcal N}(\eta(x))$:
\begin{equation}
\sum_{n=1}^{\mathcal N}\dot{\eta}(x_n){\rm v}^{(m)}_n\prod_{j\neq n}^{\mathcal N}
\left(\eta-\eta(x_j)\right)
\propto P_m(\eta),\quad m=0,1,\ldots,\mathcal{N}-1.
\label{Pmgen}
\end{equation}
\end{theo}
The  solution of the matrix equation \eqref{mateq},
$\gamma^{(m)}_n(t)= e^{i(\mathcal{E}(\mathcal{N})
-\mathcal{E}(m))t}\text{v}^{(m)}_n$
generates 
\begin{equation*}
\prod_{n=1}^{\mathcal N}\bigl(\eta(x)-\eta(x_n)\bigr)
             +\epsilon\alpha_m e^{i(\mathcal{E}(\mathcal{N})
             -\mathcal{E}(m))t} P_m(\eta(x))
             +O(\epsilon^2),
\end{equation*}
corresponding to $c_n\propto \epsilon\,\delta_{n\,m}$ 
in the polynomial solution \eqref{psiNdef}. 
Here $\alpha_m$ is a certain constant.
In other words, the present construction provides the expression of any {\em classical 
orthogonal polynomial\/} $P_m(\eta)$ in terms of the zeros of 
a higher degree polynomial
$P_{\mathcal N}(\eta)$ of the same family.

The following Lemma is well known.
\begin{lemma}\label{lemm:1}
The Lagrangian interpolation of a polynomial 
$Q(x)$ {\rm (}${\rm deg}Q=m${\rm )} by a higher degree polynomial
$\widetilde{Q}(x)$, {\rm (}${\rm deg}\widetilde{Q}=\mathcal{N}>m${\rm )} is exact:
\begin{align} 
 Q(x)&=\sum_{n=1}^{\mathcal N}\frac{Q(x_n)}{\widetilde{Q}'(x_n)}
  \cdot\left(\frac{\widetilde{Q}(x)}{x-x_n}\right),
\\
  \widetilde{Q}(x)&\eqdef d_{\mathcal N}\prod_{n=1}^{\mathcal N}\left(x-x_n\right),\quad
   \widetilde{Q}'(x_n)=d_{\mathcal N}\prod_{j\neq n}^{\mathcal N}
   \left(x_n-x_j\right).
\end{align}
\end{lemma}

In terms of the  Lemma, Theorem 3.2 can be stated as
\begin{coro}\label{coro:1}
\begin{align}
{\rm v}^{(m)}_n\propto \frac{P_m(\eta(x_n))}{\dot{\eta}(x_n)P'_{\mathcal N}(\eta(x_n))}
=\frac{P_m(\eta(x_n))}{\quad \left.\left(\frac{d P_{\mathcal N}(\eta(x))}{dx}\right)\right|_{x=x_n}},\quad n&=1,\ldots, {\mathcal N},\n[4pt]
 m&=0,1,\ldots, {\mathcal N}-1.
\label{coro1eq}
\end{align}
\end{coro}

\begin{rema}\label{rem:0}
This is rather remarkable. The matrix $\mathcal{M}$ \eqref{Mdef},
constructed by the zeros of $P_{\mathcal N}(\eta(x))$ and the basic parameters of $\widetilde{\mathcal H}$ \eqref{tevo} only, contains all the information of the values of
lower degree polynomials at these zeros $\{P_m(\eta(x_n))\}$, $m=0,1,\ldots,{\mathcal N}-1$ as eigenvectors.
\end{rema}

\begin{rema}\label{rem:1}
The above eigenvalues \eqref{eigvals} are {\rm algebraic numbers} based on the basic
parameters of $\widetilde{\mathcal H}$ \eqref{Htil}, $\mathcal{N}$ and the zeros,
 $\{\eta(x_n)\}$,  $\{\dot{\eta}(x_n)\}$. 
 The very fact that they are independent of the zeros
 means that the way that the zeros enter the matrix elements $\mathcal{M}_{n\,m}$
 is essential but that the explicit values of  $\{\eta(x_n)\}$,  
 $\{\dot{\eta}(x_n)\}$ are 
 irrelevant.
However, their explicit values are indispensable for the exact values of the eigenvectors 
$\{{\rm v}^{(m)}_n\}$ to reproduce the lower degree polynomials 
$\{P_m(\eta)\}$ \eqref{Pmgen}.
The algebraic equations satisfied by the zeros play essential roles. They are simply
obtained by evaluating  the polynomial equation \eqref{HtilN} at the zeros:
\begin{equation}
0=\left.\widetilde{\mathcal H}_{\mathcal N}
P_{\mathcal N}\left(\eta(x)\right)\right|_{x=x_n}.
\label{algeq}
\end{equation}
For the Hermite, Laguerre and Jacobi polynomials, these are well known, 
see \eqref{hermal}, 
\eqref{lagal} and \eqref{jacal}. 
The matrix elements of $\mathcal{M}$ \eqref{Mdef} are very closely related with them.
\end{rema}
\begin{rema}\label{rem:2}
The matrix $\mathcal{M}$ \eqref{Mdef} is conceptually and structurally much simpler 
than the related matrices  introduced 
for the Hermite, Laguerre and Jacobi polynomials  by Ahmed et al {\rm \cite{ahmed}} and
for the Wilson and Racah by Bihun-Calogero {\rm \cite{bihun-cal}}.
The corresponding matrices have the same eigenvalues 
{\rm(}up to an additive constant and an overall factor\/{\rm)}.
As will be shown  in \S\ref{sec:oQM}, the matrices in {\rm \cite{ahmed}} for  
the Hermite, Laguerre and Jacobi polynomials and the matrix $\mathcal{M}$ 
in Theorem \ref{theo:2} share the same eigenvectors
\eqref{hermMA},\eqref{lagMB},\eqref{jacMC}. 
Based on this fact we demonstrate  explicitly that the 
eigenvectors $\{{\rm v}^{(m)}_n\}$ of $\mathcal{M}$ have 
the above form \eqref{coro1eq} for the
Classical orthogonal polynomials. 
\end{rema}

In the subsequent sections, we provide explicit examples and data of the classical
orthogonal polynomials for which the above Theorems apply.
The Hermite, Laguerre and Jacobi polynomials in section four.
They are the main part of the eigenfunctions of exactly solvable systems in ordinary quantum mechanics.
The Wilson and Askey-Wilson polynomials and their reduced form polynomials are discussed
in section five. 
They are the main part of the eigenfunctions of exactly solvable systems in 
{\em discrete quantum mechanics with pure imaginary shifts\/} \cite{os13,os24}.
The Racah and $q$-Racah polynomials and their reduced form polynomials are examined
in section six. 
They are the main part of the eigenfunctions of exactly solvable systems in 
{\em discrete quantum mechanics with real shifts\/} \cite{os12,os24}.
The examples in sections five and six are grouped according to the sinusoidal coordinates.
In section five, the group of linear in $x$, $\eta(x)=x$, 
contains the continuous Hahn and Meixner-Pollaczek.
The group of quadratic in $x$, $\eta(x)=x^2$, 
consists of the Wilson and continuous Hahn.
The group of $\eta(x)=\cos x$ comprises of the Askey-Wilson, continuous dual $q$-Hahn,
Al-Salam-Chihara, continuous (big) $q$-Hermite, and continuous $q$-Jacobi (Laguerre).
In section six, the group of linear in $x$, $\eta(x)=x$, contains the Hahn, Krawtchouk, 
Meixner and Charlier.
The group of quadratic in $x$, $\eta(x)=x(x+d)$, consists of the Racah and  dual Hahn.
The group of linear in $q^{\pm x}$, $\eta(x)=q^{-x}-1,1-q^x$, contains the $q$-Hahn,
(quantum, affine) $q$-Krawtchouk, little $q$-Jacobi, $q$-Meixner, little $q$-Laguerre,
Al-Salam-Carlitz II, and  (alternative) $q$-Charlier. 
The group of `bilinear' in $q^{\pm x}$ consists of the $q$-Racah and dual $q$-Hahn.
For each example, we will provide the explicit form of the polynomial $\{P_n(\eta)\}$,
the sinusoidal coordinate $\eta(x)$ and the second order differential/difference 
operator $\widetilde{\mathcal H}$, so that self contained verification of the Theorems
could be made, either algebraically or numerically.
In these examples, we adopt the following short hand notation
\begin{equation}
y_n\equiv \eta(x_n),\quad n=1,\ldots, \mathcal{N},
\end{equation}
for notational simplicity, except for the linear cases, $\eta(x)=x$.

%%%%%%%%%%%%%%%%%%%%%%%%%%%%%%%%%%%%%%%%%%%%%%%%%%%%%%%%%%%%%%%
%                                                             %
%  4. Examples from ordinary quantum mechanics                                 %
%                                                             %
%%%%%%%%%%%%%%%%%%%%%%%%%%%%%%%%%%%%%%%%%%%%%%%%%%%%%%%%%%%%%%%
\section{Examples from ordinary quantum mechanics  }
\label{sec:oQM}

\subsection{Hermite}
\label{sec:her}

The system has no parameter and the various data are:
\begin{align} %
&\widetilde{\mathcal H}=-\frac{d^2}{dx^2}+2x\frac{d}{dx},\quad -\infty<x<\infty,
\quad \eta(x)=x,\quad \mathcal{E}(n)=2n, \\
&\phi_0(x)^2=e^{-x^2},\quad P_n(\eta)=H_n(\eta),\quad \text{Hermite polynomial}.
\end{align}
The matrix ${\mathcal M}$ \eqref{Mdef}  reads
\begin{align} 
   {\mathcal M}_{n\,m}&=2\delta_{n\,m}\left(\mathcal{N}+
  { \sum_{j<k}^{\mathcal N}}{}^\prime\frac1{x_n-x_j}\cdot\frac1{x_n-x_k}
   -{x_n}\sum_{j=1}^{\mathcal N}{}^\prime\frac1{x_n-x_j}\right)\n
   & \quad +2(1-\delta_{n\,m})\frac1{x_n-x_m}
   \left(\sum_{j=1, \neq m}^{\mathcal N}\!\!\!{}^\prime\frac1{x_n-x_j}-x_n\right).
\end{align}
A prime appended to a sum  indicates that the singular terms are omitted.
This matrix has a remarkable Diophantine property.
It is elementary to verify for lower $\mathcal{N}$ 
that the eigenvalues of $\mathcal{M}$ 
are all integers
\begin{equation*}
2(\mathcal{N}-m),\quad m=0,1,\ldots,\mathcal{N}-1,
\end{equation*}
for {\em arbitrary distinct complex numbers} $\{x_n\}$.

The polynomial equation for $H_{\mathcal N}$ \eqref{algeq} yields the well-known
algebraic equations among the zeros $\{x_n\}$:
\begin{equation}
\sum_{j=1}^{\mathcal N}{}^\prime\frac1{x_n-x_j}=x_n.
\label{hermal}
\end{equation}
In terms of the zeros $\{x_n\}$ of $H_{\mathcal N}$,  $H_{\mathcal N}(x_n)=0$, 
Ahmed et  al \cite{ahmed}
introduced an $\mathcal{N}\times \mathcal{N}$ matrix
\begin{equation}
A_{n\,m}=\delta_{n\,m}\sum_{j=1}^{\mathcal N}{}^\prime \frac1{(x_n-x_j)^2}
-(1-\delta_{n\,m})\frac1{(x_n-x_m)^2},
\end{equation}
having the eigenvector
\begin{equation}
v^{(m)}_n=\frac{H_m(x_n)}{H_{\mathcal{N}-1}(x_n)},\quad n=1,\ldots,\mathcal{N},
\end{equation}
corresponding to the eigenvalue $\mathcal{N}-m-1$, $m=0,1,\ldots,\mathcal{N}-1$.
It is elementary to show
\begin{equation}
\mathcal{M}=2(A+1)
\label{hermMA}
\end{equation}
by using \eqref{hermal} and another equation
\begin{equation}
\sum_{j=1}^{\mathcal N}{}^\prime\frac1{(x_n-x_j)^2}
=\frac23(\mathcal{N}-1)-\frac13x_n^2.
\end{equation}
This provides the direct derivation of Corollary \ref{coro:1} \eqref{coro1eq}, 
since $H_{\mathcal N}'(x)=2H_{\mathcal{N}-1}(x)$.

%%%%%%%%%%%%%%%%%%%%%%%%%%%%%%%%%%%
\subsection{Laguerre}
\label{sec:lag}

The system has one parameter $g>-\tfrac12$ and the various data are:
\begin{align} %
&\widetilde{\mathcal H}=-\frac{d^2}{dx^2}+2(x-\frac{g}x)\frac{d}{dx},
\quad 0<x<\infty,
\quad \eta(x)=x^2,\quad \mathcal{E}(n)=4n, \\
&\phi_0(x)^2=e^{-x^2}(x^2)^g,\quad P_n(\eta)=L_n^{(\alpha)}(\eta),\quad 
\text{Laguerre polynomial},
\quad \alpha\eqdef g-\tfrac12.
\end{align}
The matrix ${\mathcal M}$ \eqref{Mdef}  reads ($y_n=x_n^2$)
\begin{align} 
   {\mathcal M}_{n\,m}&=4\delta_{n\,m}\left(\mathcal{N}+
  2y_n\sum_{j<k}^{\mathcal N}{}^\prime\frac1{y_n-y_j}\cdot\frac1{y_n-y_k}
   -(y_n-\alpha-1)\sum_{j=1}^{\mathcal N}{}^\prime\frac1{y_n-y_j}\right)\n
   & \quad +4(1-\delta_{n\,m})\frac{x_m}{x_n}\cdot\frac1{y_n-y_m}
   \left(2y_n\!\sum_{j=1, \neq m}^{\mathcal N}\!\!\!{}^\prime
   \frac1{y_n-y_j}-(y_n-\alpha-1)\right).
   \label{Mlag}
\end{align}
This matrix has a remarkable Diophantine property.
It is elementary to verify for lower $\mathcal{N}$ that the eigenvalues 
of $\mathcal{M}$ are all integers
\begin{equation*}
4(\mathcal{N}-m),\quad m=0,1,\ldots,\mathcal{N}-1,
\end{equation*}
for {\em arbitrary distinct complex numbers} $\{y_n\}$ and $\{x_n\}$ except for $0$.

The polynomial equation for $L_{\mathcal N}^{(\alpha)}$ \eqref{algeq} 
with a change of variables yields the well-known
algebraic equations among the zeros $\{y_n\}$:
\begin{equation}
y_n\sum_{j=1}^{\mathcal N}{}^\prime\frac1{y_n-y_j}=\frac12(y_n-(\alpha+1)).
\label{lagal}
\end{equation}
In terms of the zeros $\{y_n\}$ of $L_{\mathcal N}^{(\alpha)}$,  
$L_{\mathcal N}^{(\alpha)}(y_n)=0$, Ahmed et  al \cite{ahmed}
introduced an $\mathcal{N}\times \mathcal{N}$ matrix
\begin{equation}
B_{n\,m}=\delta_{n\,m}\sum_{j=1}^{\mathcal N}{}^\prime \frac{y_j}{(y_n-y_j)^2}
-(1-\delta_{n\,m})\frac{y_m}{(y_n-y_m)^2},
\end{equation}
having the eigenvector
\begin{equation}
v^{(m)}_n=\frac{L_m^{(\alpha)}(y_n)}{L_{\mathcal{N}-1}^{(\alpha)}(y_n)},\quad 
n=1,\ldots,\mathcal{N},
\end{equation}
corresponding to the eigenvalue $\tfrac12(\mathcal{N}-m-1)$, 
$m=0,1,\ldots,\mathcal{N}-1$.
It is elementary to show
\begin{equation}
\mathcal{M}=4\mathcal{D}(2B+1)\mathcal{D}^{-1},
\quad \mathcal{D}\eqdef\text{diag}(x_1,x_2,\ldots,x_{\mathcal N}),
\label{lagMB}
\end{equation}
by using \eqref{lagal} and another equation \cite{ahmed} 
\begin{align} 
y_n^2\sum_{j=1}^{\mathcal N}{}^\prime\frac1{(y_n-y_j)^2}
&=-\frac1{12}\left((\alpha+1)(\alpha+5)-2(2\mathcal{N}+\alpha+1)y_n+y_n^2\right).
\end{align}
This means that the eigenvectors of $\mathcal{M}$ \eqref{Mlag} are
\begin{equation}
\text{v}^{(m)}_n
=x_n\frac{L_m^{(\alpha)}(y_n)}{L_{\mathcal{N}-1}^{(\alpha)}(y_n)}\propto
\frac1{x_n}\frac{L_m^{(\alpha)}(y_n)}{L_{\mathcal{N}}^{(\alpha)}{}^\prime(y_n)},
\qquad m=0,1,\ldots,\mathcal{N}-1,
\end{equation}
providing the direct derivation of Corollary \ref{coro:1} \eqref{coro1eq}. 
In deriving  the 
final proportionality relation, the following identity 
((5.1.14) of \cite{szego}) is useful:
\begin{align}
&\eta\frac{d L^{(\alpha)}_n(\eta)}{d\eta}
=-\eta L_{n-1}^{(\alpha+1)}(\eta)
=nL_{n}^{(\alpha)}(\eta)-(n+\alpha)L_{n-1}^{(\alpha)}(\eta).
\end{align}

%%%%%%%%%%%%%%%%%%%%%%%%%%%%%%%%%%%
\subsection{Jacobi}
\label{sec:jac}

The system has two parameters $g>-\tfrac12$, $h>-\tfrac12$ 
and the various data are:
\begin{align} %
&\widetilde{\mathcal H}=-\frac{d^2}{dx^2}-2(g\cot x-h\tan x)\frac{d}{dx},\quad 
0<x<\frac{\pi}2,
\\
& \quad \eta(x)=\cos2x,\quad \dot{\eta}(x)=-2\sin2x, \qquad 
(\dot{\eta}(x))^2=4(1-\eta(x)^2),\\  
&\phi_0(x)^2=(\sin^2x)^g(\cos^2x)^h,\quad P_n(\eta)
=P_n^{(\alpha,\beta)}(\eta),\quad \text{Jacobi polynomial},\\
&\quad \mathcal{E}(n)=4n(n+g+h)=4n(n+\alpha+\beta+1),\quad 
\alpha\eqdef g-\tfrac12, \ \beta\eqdef h-\tfrac12.
\end{align}
The matrix ${\mathcal M}$ \eqref{Mdef}  reads ($y_n=\cos2x_n$)
\begin{align} 
   {\mathcal M}_{n\,m}&=4\delta_{n\,m}\left(\mathcal{N}(\mathcal{N}+\alpha+\beta+1)+
  2(1-y_n^2)\sum_{j<k}^{\mathcal N}{}^\prime\frac1{y_n-y_j}
  \cdot\frac1{y_n-y_k}\right.\n
  & \hspace{70mm}\left. -((\alpha+\beta)y_n
  +\alpha-\beta)\sum_{j=1}^{\mathcal N}{}^\prime\frac1{y_n-y_j}\right)\n
   & \quad +4(1-\delta_{n\,m})\frac{\sin2x_m}{\sin2x_n}\cdot\frac1{y_n-y_m}
   \left(2(1-y_n^2)\!\sum_{j=1, \neq m}^{\mathcal N}\!\!\!{}^\prime\frac1{y_n-y_j}
   -((\alpha+\beta)y_n+\alpha-\beta)\right).
   \label{Mjac}
\end{align}
This matrix has a remarkable Diophantine property when $\alpha+\beta$ is an integer.
It is elementary to verify for lower $\mathcal{N}$ that the eigenvalues of 
$\mathcal{M}$ \eqref{Mjac} are
\begin{equation*}
4(\mathcal{N}-m)(\mathcal{N}+m+\alpha+\beta+1),\quad m=0,1,\ldots,\mathcal{N}-1,
\end{equation*}
for {\em arbitrary distinct complex numbers} 
$\{y_n\}$ and $\{x_n\}$ except for $0$ mod $\pi/2$. 

The polynomial equation for $P_{\mathcal N}^{(\alpha,\beta)}$ \eqref{algeq} 
with a change of variables yields the well-known
algebraic equations among the zeros $\{y_n\}$:
\begin{equation}
(1-y_n^2)\sum_{j=1}^{\mathcal N}{}^\prime\frac1{y_n-y_j}
=\frac12\bigl((\alpha+\beta+2)y_n+\alpha-\beta\bigr).
\label{jacal}
\end{equation}
In terms of the zeros $\{y_n\}$ of $P_{\mathcal N}^{(\alpha,\beta)}$,  
$P_{\mathcal N}^{(\alpha,\beta)}(y_n)=0$, Ahmed et  al \cite{ahmed}
introduced an $\mathcal{N}\times \mathcal{N}$ matrix
\begin{equation}
C_{n\,m}=\delta_{n\,m}\sum_{j=1}^{\mathcal N}{}^\prime 
\frac{(1-y_j^2)}{(y_n-y_j)^2}
-(1-\delta_{n\,m})\frac{(1-y_m^2)}{(y_n-y_m)^2},
\end{equation}
having the eigenvector
\begin{equation}
v^{(m)}_n=\frac{P_m^{(\alpha,\beta)}(y_n)}{P_{\mathcal{N}-1}^{(\alpha,\beta)}(y_n)},\quad 
n=1,\ldots,\mathcal{N},
\end{equation}
corresponding to the eigenvalue 
$\tfrac12(\mathcal{N}-m-1)(\mathcal{N}+m+\alpha+\beta)$, 
$m=0,1,\ldots,\mathcal{N}-1$.
It is elementary to show
\begin{equation}
\mathcal{M}=4\mathcal{D}(2C+2\mathcal{N}+\alpha+\beta)\mathcal{D}^{-1},\quad 
\mathcal{D}\eqdef\text{diag}(\sin2x_1,\sin2x_2,\ldots,\sin2x_{\mathcal N}),
\label{jacMC}
\end{equation}
by using \eqref{jacal} and another equation \cite{ahmed} 
\begin{align} 
&(1-y_n^2)^2\sum_{j=1}^{\mathcal N}{}^\prime\frac1{(y_n-y_j)^2}\n
&\qquad =\frac1{3}(\mathcal{N}-1)
(\mathcal{N}+\alpha+\beta+2)-\frac1{12}(\alpha-\beta)^2
-\frac16(\alpha-\beta)(\alpha+\beta+6)y_n\n
&\qquad \quad -\frac1{12}\left[4\mathcal{N}(\mathcal{N}+\alpha+\beta+1)
+(\alpha+\beta+2)(\alpha+\beta+6)\right]y_n^2.
\end{align}
This means that the eigenvectors of $\mathcal{M}$ \eqref{Mjac} are
\begin{equation}
\text{v}^{(m)}_n=\sin2x_n\frac{P_m^{(\alpha,\beta)}(y_n)}
{P_{\mathcal{N}-1}^{(\alpha,\beta)}(y_n)}\propto
\frac1{\sin2x_n}\frac{P_m^{(\alpha,\beta)}(y_n)}
{P_{\mathcal{N}}^{(\alpha,\beta)}{}^\prime(y_n)},
\qquad m=0,1,\ldots,\mathcal{N}-1,
\end{equation}
providing the direct derivation of Corollary \ref{coro:1} \eqref{coro1eq}. 
In deriving  the 
final proportionality relation, the following identity 
((4.5.7) of \cite{szego}) is useful:
\begin{align}
&(2n+\alpha+\beta)(1-\eta^2)\frac{d P^{(\alpha,\beta)}_n(\eta)}{d\eta}\n
&=-n\left[(2n+\alpha+\beta)\eta+\beta-\alpha\right]
P_{n}^{(\alpha,\beta)}(\eta)+2(n+\alpha)(n+\beta)P_{n-1}^{(\alpha,\beta)}(\eta).
\end{align}

%%%%%%%%%%%%%%%%%%%%%%%%%%%%%%%%%%%%%%%%%%%%%%%%%%%%%%%%%%%%%%%
%                                                             %
%  5. Examples from discrete quantum mechanics with pure imaginary shifts                             %
%                                                             %
%%%%%%%%%%%%%%%%%%%%%%%%%%%%%%%%%%%%%%%%%%%%%%%%%%%%%%%%%%%%%%%
\section{Examples from discrete quantum mechanics with pure imaginary shifts}
\label{sec:dQMpi}

The difference operator $\widetilde{\mathcal{H}}$ 
governing the classical orthogonal polynomials 
\eqref{HtilPn}
belonging  to this class
 depends on an analytic function $V(x)$ of $x$:
\begin{align}
  \widetilde{\mathcal{H}}  
  &=V(x)(e^{-i\gamma \partial_x}-1)+V^*(x)(e^{i\gamma \partial_x}-1).
  \label{HVV*}
\end{align}
Throughout this section we use $i\eqdef\sqrt{-1}$.
The function $V^*(x)$ 
is an analytic function of $x$ obtained from $V(x)$ by the $*$-operation,
which is defined as follows.
If $f(x)=\sum\limits_{n}a_nx^n$, $a_n\in\mathbb{C}$, then
$f^*(x)\eqdef\sum\limits_{n}a_n^*x^n$, in which $a_n^*$ is the complex
conjugation of $a_n$. Here $\gamma$ is a real parameter specifying 
the shifts of the functions:
\begin{equation*}
e^{\pm i\gamma\partial_x}\psi(x)=\psi(x\pm i\gamma).
\end{equation*}
We consider two different cases $\gamma=1$ for the Wilson polynomials 
and its reduced form polynomials
and $\gamma=\log q$ and $0<q<1$ for the Askey-Wilson polynomials 
and its reduced form polynomials.
The polynomials are assembled into three groups according 
to the form of the sinusoidal coordinate 
$\eta(x)=x,x^2$ and $\cos x$. In each group, we start from 
the most generic member and move to simpler ones.
For a comprehensive exposition of these polynomials 
in the discrete quantum mechanics formulation, 
see \cite{os13,os24}.
%%%%%%%%%%%%%%%%
\subsection{Polynomials having  $\eta(x)=x$, $-\infty<x<\infty$, $\gamma=1$}
\label{sec:lineta}
This group consists of two members, the continuous Hahn \S\ref{sec:contH} 
and the Meixner-Pollaczek
\S\ref{sec:Meix-Pol}.
The matrix $\mathcal{M}$ \eqref{Mdef} for this group reads
\begin{align} 
\mathcal{M}_{n\,m}&=\delta_{n\,m}\left(
 -\frac{V(x_n)\prod_{j\neq n}^{\mathcal N}(x_n-i-x_j)+
V^*(x_n)\prod_{j\neq n}^{\mathcal N}(x_n+i-x_j)}
{\prod_{j\neq n}^{\mathcal N}(x_n-x_j)}\right.\n[4pt]
&\hspace{70mm}
\left. +\mathcal{E}(\mathcal{N})+V(x_n)+V^*(x_n)
\phantom{\frac{(x_n)}{x_n}}\hspace*{-9mm}\right)\n
&\quad +(1-\delta_{n\,m})\frac{i\left(V(x_n)\prod_{j\neq n,m}^{\mathcal N}(x_n-i-x_j)-
V^*(x_n)\prod_{j\neq n,m}^{\mathcal N}(x_n+i-x_j)\right)}
{\prod_{j\neq n}^{\mathcal N}(x_n-x_j)}.
\label{McontH}
\end{align}
It has eigenvalues $\mathcal{E}(\mathcal{N})-\mathcal{E}(m)$, 
$m=0,1,\ldots,\mathcal{N}-1$,
for arbitrary distinct complex values of $\{x_n\}$.
For the zeros $\{x_n\}$ of $P_{\mathcal N}$, $P_{\mathcal N}(x_n)=0$, 
it is  straightforward
to verify Theorem \ref{theo:2} numerically for lower $\mathcal{N}$. 
That is,  the eigenvectors of the matrix 
$\mathcal{M}$ \eqref{McontH} generate
the lower degree polynomials $\{P_m(x)\}$, $m=0,1,\ldots,\mathcal{N}-1$ 
as in \eqref{Pmgen}.
%%%%%%%%%%%%%%%%%%%%%
\subsubsection{continuous Hahn}
\label{sec:contH}
This polynomial depends on two complex parameters 
$(a_1,a_2)$ and the various data are \cite{koekswart,os13}:
\begin{align} 
V(x)&=(a_1+ix)(a_2+ix),\quad V^*(x)=(a_3-ix)(a_4-ix),\\
\mathcal{E}(n)&=n(n+b_1-1),\quad b_1\eqdef\sum_{j=1}^4a_j,
\quad  \{a_3,a_4\}=\{a_1^*,a_2^*\}\ \text{as a set},
\quad \text{Re}\,a_j>0,\\
\phi_0(x)^2&=\prod_{j=1}^2\Gamma(a_j+ix)\Gamma(a_j^*-ix),\\
P_n(x)&=
  i^n\frac{(a_1+a_3)_n(a_1+a_4)_n}{n!}\,
  {}_3F_2\Bigl(\genfrac{}{}{0pt}{}{-n,\,n+b_1-1,\,a_1+ix}
  {a_1+a_3,\,a_1+a_4}\!\Bigm|\!1\Bigr).
 \label{defcH}
\end{align}
The polynomial equation for $P_{\mathcal N}$ \eqref{algeq} 
provides algebraic equations for the zeros $\{x_n\}$:
\begin{align} 
(a_1+ix_n)(a_2+ix_n)\prod_{j\neq n}^\mathcal{N}(x_n-i-x_j)=
(a_3-ix_n)(a_4-ix_n)\prod_{j\neq n}^\mathcal{N}(x_n+i-x_j).
\end{align}
%%%%%%%%%%%%%%%%%%%%%
\subsubsection{Meixner-Pollaczek}
\label{sec:Meix-Pol}
This polynomial depends on two real parameters $(a,\phi)$ 
and the various data are \cite{koekswart,os13}:
\begin{align} 
V(x)&\eqdef e^{i(\tfrac{\pi}2-\phi)}(a+ix),
\quad V^*(x)=e^{-i(\tfrac{\pi}2-\phi)}(a-ix),\\
\mathcal{E}(n)&=2n\sin\phi,\quad 
\quad a>0,\quad 0<\phi<\pi\\
\phi_0(x)^2&=e^{(2\phi-\pi)x}\Gamma(a+ix)\Gamma(a-ix),\\
 P_n(x)&=\frac{(2a)_n}{n!}\,e^{in\phi}
  {}_2F_1\Bigl(\genfrac{}{}{0pt}{}{-n,\,a+ix}{2a}\Bigm|
  1-e^{-2i\phi}\Bigr).
  \label{defMP}
\end{align}
The polynomial equation for $P_{\mathcal N}$ \eqref{algeq} 
provides algebraic equations for the zeros $\{x_n\}$:
\begin{align} 
e^{i(\tfrac{\pi}2-\phi)x}(a+ix_n)\prod_{j\neq n}^\mathcal{N}(x_n-i-x_j)=
e^{-i(\tfrac{\pi}2-\phi)x)
}(a-ix_n)\prod_{j\neq n}^\mathcal{N}(x_n+i-x_j),
\end{align}
which simplifies for $\phi=\tfrac{\pi}2$.

%%%%%%%%%%%%%%%%
\subsection{Polynomials having  $\eta(x)=x^2$, $0<x<\infty$, $\gamma=1$}
\label{sec:quadeta}
This group consists of two members, the Wilson \S\ref{sec:Wilson} 
and the continuous dual Hahn
\S\ref{sec:contdH}.
The matrix $\mathcal{M}$ \eqref{Mdef} for this group reads
\begin{align} 
\mathcal{M}_{n\,m}&=\delta_{n\,m}\left(
 -\frac{V(x_n)\prod_{j\neq n}^{\mathcal N}\bigl((x_n-i)^2-y_j\bigr)+
V^*(x_n)\prod_{j\neq n}^{\mathcal N}\bigl((x_n+i)^2-y_j\bigr)}
{\prod_{j\neq n}^{\mathcal N}(y_n-y_j)}\right.\n[4pt]
&\hspace{70mm}
\left. +\mathcal{E}(\mathcal{N})+V(x_n)+V^*(x_n)
\phantom{\frac{(x_n)}{x_n}}\hspace*{-9mm}\right)\n
&\quad +(1-\delta_{n\,m})\frac{x_m}{x_n}
\frac1{\prod_{j\neq n}^{\mathcal N}(y_n-y_j)}
\left(V(x_n)(1+2ix_n)\prod_{j\neq n,m}^{\mathcal N}\bigl((x_n-i)^2-y_j\bigr)\right.\n
&  \hspace{60mm} \left. +
V^*(x_n)(1-2ix_n)\prod_{j\neq n,m}^{\mathcal N}\bigl((x_n+i)^2-y_j\bigr)\right).
\label{MWil}
\end{align}
It has eigenvalues $\mathcal{E}(\mathcal{N})-\mathcal{E}(m)$, 
$m=0,1,\ldots,\mathcal{N}-1$,
for arbitrary distinct complex values of $\{x_n\}$  
except for the poles of $V$ and $V^*$ and $\{y_n=x_n^2\}$. 
They are all integers for the continuous dual Hahn.
The same for the Wilson, if $b_1$ \eqref{ctdHphi0} is an integer.
For the zeros $\{y_n=x_n^2\}$ of $P_{\mathcal N}$, $P_{\mathcal N}(y_n)=0$, 
it is  straightforward
to verify Theorem \ref{theo:2} numerically for lower $\mathcal{N}$. 
That is,  the eigenvectors of the matrix 
$\mathcal{M}$ \eqref{MWil} generate
the lower degree polynomials $\{P_m(\eta)\}$, $m=0,1,\ldots,\mathcal{N}-1$ 
as in \eqref{Pmgen}.

%%%%%%%%%%%%%%%%%%%%
\subsubsection{Wilson}
\label{sec:Wilson}
This polynomial depends on four real parameters or two complex parameters, 
$ \{a_1^*,a_2^*,a_3^*,a_4^*\}=\{a_1,a_2,a_3,a_4\}$ as a set, 
and the various data are \cite{koekswart,os13}:
\begin{align} 
   V(x)&=
  \frac{(a_1+ix)(a_2+ix)(a_3+ix)(a_4+ix)}{2ix(2ix+1)}, 
  \quad V^*(x)=V(-x),\quad
 \text{Re}\,a_j>0, \\ 
  \phi_0(x)^2&=
 \frac{\prod_{j=1}^4\Gamma(a_j+ix)\Gamma(a_j-ix)}{\Gamma(2ix)\Gamma(-2ix)},
 \quad  \mathcal{E}(n)=n(n+b_1-1),\quad  b_1\eqdef\sum_{j=1}^4a_j,
  \label{Wilsonphi0}\\[4pt]
  P_n(\eta(x))
  &=(a_1+a_2)_n(a_1+a_3)_n(a_1+a_4)_n
  \times {}_4F_3\Bigl(
  \genfrac{}{}{0pt}{}{-n,\,n+b_1-1,\,a_1+ix,\,a_1-ix}
  {a_1+a_2,\,a_1+a_3,\,a_1+a_4}\Bigm|1\Bigr).
  \label{defW}
\end{align}
The polynomial equation for $P_{\mathcal N}$  \eqref{algeq} provides 
algebraic equations for the zeros $\{y_n\}$:
\begin{align} 
  \prod_{k=1}^4(a_k+ix_n)\cdot\prod_{j\neq n}^{\mathcal N}\bigl((x_n-i)^2-y_j\bigr)
  =\prod_{k=1}^4(a_k-ix_n)\cdot\prod_{j\neq n}^{\mathcal N}\bigl((x_n+i)^2-y_j\bigr).
\end{align}
All the non-$q$ polynomials in this section can be obtained from the Wilson by reductions.

\subsubsection{continuous dual Hahn}
\label{sec:contdH}
This is a restricted case of the Wilson polynomial with $a_4=0$.
The parameters are restricted by $\{a_1^*,a_2^*,a_3^*\}=\{a_1,a_2,a_3\}$, as a set 
and  the various data are \cite{koekswart,os13}:
\begin{align} 
   V(x)&=
  \frac{(a_1+ix)(a_2+ix)(a_3+ix)}{2ix(2ix+1)}, \quad V^*(x)=V(-x),\quad
 \text{Re}\,a_j>0, \\ 
  \phi_0(x)^2&=
 \frac{\prod_{j=1}^3\Gamma(a_j+ix)\Gamma(a_j-ix)}{\Gamma(2ix)\Gamma(-2ix)},
 \quad  \mathcal{E}(n)=n,
  \label{ctdHphi0}\\[4pt]
  P_n(\eta(x))
  &=(a_1+a_2)_n(a_1+a_3)_n  \times {}_3F_2\Bigl(
  \genfrac{}{}{0pt}{}{-n,,\,a_1+ix,\,a_1-ix}
  {a_1+a_2,\,a_1+a_3}\Bigm|1\Bigr).
  \label{defcontdH}
\end{align}
The polynomial equation for $P_{\mathcal N}$ \eqref{algeq} provides 
algebraic equations for the zeros $\{y_n\}$:
\begin{align} 
  \prod_{k=1}^3(a_k+ix_n)\cdot\prod_{j\neq n}^{\mathcal N}\bigl((x_n-i)^2-y_j\bigr)
  =\prod_{k=1}^4(a_k-ix_n)\cdot\prod_{j\neq n}^{\mathcal N}\bigl((x_n+i)^2-y_j\bigr).
\end{align}

%%%%%%%%%%%%%%%%
\subsection{Polynomials having  $\eta(x)=\cos x$, $0<x<\pi$, $e^\gamma=q$}
\label{sec:cos}
The Askey-Wilson polynomial and its six reduced case polynomials belong to this group.
For this group, let us introduce new symbols related 
with  the zeros $\{x_n\}$, ($\dot{\eta}(x)=-\sin x$):
\begin{equation}
z_n\eqdef e^{ix_n}=\cos x_n+i\sin x_n=y_n-i\dot{\eta}(x_n),
\quad n=1,\ldots, \mathcal{N}.
\label{zndef}
\end{equation}
The matrix $\mathcal{M}$ \eqref{Mdef} for this group reads
\begin{align} 
\mathcal{M}_{n\,m}&=\delta_{n\,m}\left\{
\phantom{\frac{(x_n)}{x_n}}\hspace*{-9mm}
\mathcal{E}(\mathcal{N})+V(x_n)+V^*(x_n)\right. \n
& \qquad\qquad 
-\frac1{\prod_{j\neq n}^{\mathcal N}(y_n-y_j)}
\left(V(x_n)\prod_{j\neq n}^{\mathcal N}
\bigl((qz_n+q^{-1}z_n^{-1})/2-y_j\bigr)\right.\n
&\hspace{50mm}\left. \left. +V^*(x_n)\prod_{j\neq n}^{\mathcal N}
\bigl((q^{-1}z_n+qz_n^{-1})/2-y_j\bigr)\right)\right\}\n
&\quad +(1-\delta_{n\,m})\frac{\sin2x_m}{\sin2x_n}
\frac{(q^{-1}-1)}{2\prod_{j\neq n}^{\mathcal N}(y_n-y_j)}\n
&
\qquad \times \left(V(x_n)z_n^{-1}(1-qz_n^{2})
\prod_{j\neq n,m}^{\mathcal N}\bigl((qz_n+q^{-1}z_n^{-1})/2-y_j\bigr)\right.\n
&  \hspace{15mm} \left. +
V^*(x_n)z_n(1-qz_n^{-2})\prod_{j\neq n,m}^{\mathcal N}
\bigl((q^{-1}z_n+qz_n^{-1})/2-y_j\bigr)\right).
\label{MAW} 
\end{align}
It has eigenvalues $\mathcal{E}(\mathcal{N})-\mathcal{E}(m)$, 
$m=0,1,\ldots,\mathcal{N}-1$,
for arbitrary distinct complex values of $\{x_n\}$ 
except for the poles of $V$ and $V^*$
 and $\{y_n=\cos x_n\}$. 
For the zeros $\{y_n=\cos x_n\}$ of $P_{\mathcal N}$, $P_{\mathcal N}(y_n)=0$, 
it is  straightforward
to verify Theorem \ref{theo:2} numerically for lower $\mathcal{N}$. 
That is,  the eigenvectors of the matrix 
$\mathcal{M}$ \eqref{MAW} generate
the lower degree polynomials $\{P_m(\eta)\}$, $m=0,1,\ldots,\mathcal{N}-1$ 
as in \eqref{Pmgen}.

%%%%%%%%%%%%%%%%%%%%
\subsubsection{Askey-Wilson}
\label{sec:As-Wil}

The Askey-Wilson polynomial is the most general one with the maximal
number of four real parameters, or two complex parameters,
$ \{a_1^*,a_2^*,a_3^*,a_4^*\}=\{a_1,a_2,a_3,a_4\}$ as a set, on top of $q$.
All the other polynomials in this  group are obtained
by restricting the parameters $a_1$,\ldots,$a_4$, in one way or another.
The various data are \cite{koekswart,os13}:
\begin{align} 
   V(x)&=
  \frac{(1-a_1e^{ix})(1-a_2e^{ix})(1-a_3e^{ix})(1-a_4e^{ix})}
  {(1-e^{2ix})(1-qe^{2ix})}, \quad V^*(x)=V(-x),\quad
 |a_j|<1, \\ 
  \phi_0(x)^2&=
\frac{(e^{2ix}\,;q)_{\infty}(e^{-2ix}\,;q)_{\infty}}
  {\prod_{j=1}^4(a_je^{ix}\,;q)_{\infty}(a_je^{-ix}\,;q)_{\infty}},
  \quad \mathcal{E}(n)=(q^{-n}-1)(1-b_4q^{n-1}),
  \label{AWilsonphi0}\\
  & b_4\eqdef a_1a_2a_3a_4,\\
    P_n(\eta(x))&= a_1^{-n}(a_1a_2,a_1a_3,a_1a_4\,;q)_n
  \times
  {}_4\phi_3\Bigl(\genfrac{}{}{0pt}{}{q^{-n},\,b_4q^{n-1},\,
  a_1e^{ix},\,a_1e^{-ix}}{a_1a_2,\,a_1a_3,\,a_1a_4}\!\!\Bigm|\!q\,;q\Bigr).
  \label{defAW}
\end{align}
The polynomial equation for $P_{\mathcal N}$ \eqref{algeq} provides 
algebraic equations for the zeros $\{y_n\}$:
\begin{align} 
&\prod_{k=1}^4(1-a_kz_n)\cdot 
z_n^{-2}\prod_{j\neq n}^{\mathcal N}\bigl((qz_n+q^{-1}z_n^{-1})/2-y_j\bigr)\n
&=\prod_{k=1}^4(1-a_kz_n^{-1})\cdot z_n^{2}
\prod_{j\neq n}^{\mathcal N}\bigl((q^{-1}z_n+qz_n^{-1})/2-y_j\bigr).
\end{align}
%%%%%%%%%%%%%%%%%%%%
\subsubsection{Continuous dual $q$-Hahn}
\label{sec:contdqH}
The continuous dual $q$-Hahn polynomial is obtained by restricting
$a_4=0$ in the Askey-Wilson polynomial \S\ref{sec:As-Wil}.
The various data are \cite{koekswart,os13}: 
\begin{align} 
   V(x)&=
  \frac{(1-a_1e^{ix})(1-a_2e^{ix})(1-a_3e^{ix})}{(1-e^{2ix})(1-qe^{2ix})},
   \quad V^*(x)=V(-x),\quad
 |a_j|<1, \\ 
  \phi_0(x)^2&=
\frac{(e^{2ix}\,;q)_{\infty}(e^{-2ix}\,;q)_{\infty}}
  {\prod_{j=1}^3(a_je^{ix}\,;q)_{\infty}(a_je^{-ix}\,;q)_{\infty}},
  \quad \mathcal{E}(n)=q^{-n}-1,
  \label{cdqHphi0}\\
  &  \{a_1^*,a_2^*,a_3^*\}=\{a_1,a_2,a_3\} \quad (\text{as a set}),
  \\
    P_n(\eta(x))&= a_1^{-n}(a_1a_2,a_1a_3\,;q)_n
  \times
  {}_3\phi_2\Bigl(\genfrac{}{}{0pt}{}{q^{-n},\,
  a_1e^{ix},\,a_1e^{-ix}}{a_1a_2,\,a_1a_3}\!\!\Bigm|\!q\,;q\Bigr).
  \label{defcontdqH}
\end{align}

%%%%%%%%%%%%%%%%%%%%
\subsubsection{Al-Salam-Chihara}
\label{sec:ASC}
This is a further restriction of the continuous dual $q$-Hahn
polynomial \S\ref{sec:contdqH} with $a_3=0$.
The various data are \cite{koekswart,os13}: 
\begin{align} 
   V(x)&=
  \frac{(1-a_1e^{ix})(1-a_2e^{ix})}{(1-e^{2ix})(1-qe^{2ix})}, 
  \quad V^*(x)=V(-x),\quad
 |a_j|<1, \\ 
  \phi_0(x)^2&=
\frac{(e^{2ix}\,;q)_{\infty}(e^{-2ix}\,;q)_{\infty}}
  {\prod_{j=1}^2(a_je^{ix}\,;q)_{\infty}(a_je^{-ix}\,;q)_{\infty}},
  \quad \mathcal{E}(n)=q^{-n}-1,
  \label{ASCphi0}\\
  &  \{a_1^*,a_2^*\}=\{a_1,a_2\} \quad (\text{as a set}),
  \\
    P_n(\eta(x))&= a_1^{-n}(a_1a_2\,;q)_n
  \times
  {}_3\phi_2\Bigl(\genfrac{}{}{0pt}{}{q^{-n},\,
  a_1e^{ix},\,a_1e^{-ix}}{a_1a_2,\,0}\!\!\Bigm|\!q\,;q\Bigr).
  \label{defASC}
\end{align}

%%%%%%%%%%%%%%%%%%%%
\subsubsection{Continuous big $q$-Hermite}
\label{sec:cbqH}
This is a further restriction of the Al-Salam-Chihara polynomial
\S\ref{sec:ASC} with $a_2=0$ and it depends on one real parameter $a$ on top of $q$.
The various data are \cite{koekswart,os13}: 
\begin{align} 
   V(x)&=
  \frac{(1-a\,e^{ix})}{(1-e^{2ix})(1-qe^{2ix})}, \quad V^*(x)=V(-x),\quad
-1<a<1, \\ 
  \phi_0(x)^2&=
\frac{(e^{2ix}\,;q)_{\infty}(e^{-2ix}\,;q)_{\infty}}
  {(a\,e^{ix}\,;q)_{\infty}(a\,e^{-ix}\,;q)_{\infty}},
  \quad \mathcal{E}(n)=q^{-n}-1,
  \label{cbqHphi0}\\
    P_n(\eta(x))&= a^{-n}
  \times
  {}_3\phi_2\Bigl(\genfrac{}{}{0pt}{}{q^{-n},\,
  ae^{ix},\,ae^{-ix}}{0,\,0}\!\!\Bigm|\!q\,;q\Bigr).
  \label{defcbqH}
\end{align}

%%%%%%%%%%%%%%%%%%%%
\subsubsection{Continuous $q$-Hermite}
\label{sec:cqH}
This is a $q$ analogue of the Hermite polynomial, depending on $q$ only.
It also provides the simplest dynamical realisation of the 
$q$-oscillator algebra \cite{os11}.
The various data are \cite{koekswart,os13}: 
\begin{align} 
   V(x)&=
  \frac1{(1-e^{2ix})(1-qe^{2ix})}, \quad V^*(x)=V(-x), \\ 
  \phi_0(x)^2&=
(e^{2ix}\,;q)_{\infty}(e^{-2ix}\,;q)_{\infty},\quad \mathcal{E}(n)=q^{-n}-1,
  \label{cqHphi0}\\
    P_n(\eta(x))&= e^{inx}\,
  {}_2\phi_0\Bigl(\genfrac{}{}{0pt}{}{q^{-n},\,0}{-}
  \Bigm|q\,;q^ne^{-2ix}\Bigr).
  \label{defcqHe}
\end{align}
The polynomial equation for $P_{\mathcal N}$ \eqref{algeq} provides 
simple algebraic equations for the zeros $\{y_n\}$:
\begin{align} 
& z_n^{-2}\prod_{j\neq n}^{\mathcal N}\bigl((qz_n+q^{-1}z_n^{-1})/2-y_j\bigr)
= z_n^{2}\prod_{j\neq n}^{\mathcal N}\bigl((q^{-1}z_n+qz_n^{-1})/2-y_j\bigr).
\end{align}

%%%%%%%%%%%%%%%%%%%%
\subsubsection{Continuous $q$-Jacobi}
\label{sec:cqJ}
This polynomial depends on two real parameters $\alpha$ 
and $\beta$ and other data are \cite{koekswart,os13}: 
\begin{align} 
V(x)&=
  \frac{(1-q^{\frac12(\alpha+\frac12)}e^{ix})(1-q^{\frac12(\alpha+\frac32)}e^{ix})
  (1+q^{\frac12(\beta+\frac12)}e^{ix})(1+q^{\frac12(\beta+\frac32)}e^{ix})}
  {(1-e^{2ix})(1-qe^{2ix})},\\
 & V^*(x)=V(-x),\quad  \alpha,\beta\ge -\frac{1}{2},
 \quad  \mathcal{E}(n)=(q^{-n}-1)(1-q^{n+\alpha+\beta+1}),\\
 \phi_0(x)^2&=
\frac{(e^{2ix}\,;q)_{\infty}(e^{-2ix}\,;q)_{\infty}}
  {(q^{\frac12(\alpha+\frac12)}e^{ix},
   -q^{\frac12(\beta+\frac12)}e^{ix}\,;q^{\frac12})_{\infty}
   (q^{\frac12(\alpha+\frac12)}e^{-ix},
   -q^{\frac12(\beta+\frac12)}e^{-ix}\,;q^{\frac12})_{\infty}},
  \label{qJphi0}\\[2pt]
   P_n(\eta(x))&= 
     \frac{(q^{\alpha+1}\,;q)_n}{(q\,;q)_n}\,
  {}_4\phi_3\Bigl(\genfrac{}{}{0pt}{}{q^{-n},\,q^{n+\alpha+\beta+1},\,
  q^{\frac12(\alpha+\frac12)}e^{ix},\,q^{\frac12(\alpha+\frac12)}e^{-ix}}
  {q^{\alpha+1},\,-q^{\frac12(\alpha+\beta+1)},\,
  -q^{\frac12(\alpha+\beta+2)}}\Bigm|q\,;q\Bigr).
 \label{defcqJ}
  \end{align}
  
  %%%%%%%%%%%%%%%%%%%%
\subsubsection{Continuous $q$-Laguerre}
\label{sec:cqL}
This is a further restriction ($\beta\to\infty$ or $q^\beta\to0$) of 
the continuous $q$-Jacobi polynomial \S\ref{sec:cqJ}.
Many formulas are  simplified \cite{koekswart,os13}: 
\begin{align} 
V(x)&=
  \frac{(1-q^{\frac12(\alpha+\frac12)}e^{ix})(1-q^{\frac12(\alpha+\frac32)}e^{ix})}
{(1-e^{2ix})(1-qe^{2ix})},\quad
  V^*(x)=V(-x),\quad  \alpha\ge -\frac{1}{2},\\
 \phi_0(x)^2&=
\frac{(e^{2ix}\,;q)_{\infty}(e^{-2ix}\,;q)_{\infty}}
  {(q^{\frac12(\alpha+\frac12)}e^{ix}\,;q^{\frac12})_{\infty}
   (q^{\frac12(\alpha+\frac12)}e^{-ix}\,;q^{\frac12})_{\infty}},
   \quad  \mathcal{E}(n)=q^{-n}-1,
  \label{qLphi0}\\[2pt]
   P_n(\eta(x))&= 
 \frac{(q^{\alpha+1}\,;q)_n}{(q\,;q)_n}\,
  {}_3\phi_2\Bigl(\genfrac{}{}{0pt}{}{q^{-n},\,
  q^{\frac12(\alpha+\frac12)}e^{ix},\,q^{\frac12(\alpha+\frac12)}e^{-ix}}
  {q^{\alpha+1},\,0}\Bigm|q\,;q\Bigr).
 \label{defcqL}
 \end{align}
 %%%%%%%%%%%%%%%%%%%%%%%%%%%%%%%%%%%%%%%%%%%%%%%%%%%%%%%%%%%%%%%
%                                                             %
%  6. Examples from discrete quantum mechanics with real shifts                             %
%                                                             %
%%%%%%%%%%%%%%%%%%%%%%%%%%%%%%%%%%%%%%%%%%%%%%%%%%%%%%%%%%%%%%%
\section{Examples from discrete quantum mechanics with real shifts}
\label{sec:dQMre}

The  polynomials belonging to this class are also called 
{\em classical orthogonal polynomials of a discrete variable\/} \cite{nikiforov}.
The second order difference operator $\widetilde{\mathcal{H}}$ 
governing these polynomials 
\eqref{HtilPn} reads
\begin{align}
  \widetilde{\mathcal{H}}&=B(x)(1-e^{\partial})+D(x)(1-e^{-\partial}).
  \label{realHt}
\end{align}
Here  two non-negative functions $B(x)\ge0, D(x)\ge0$ of $x$ 
 are defined on non-negative integer 
lattice points, finite
 $[0,1,\ldots, N]$, or infinite $[0,1,\ldots,\infty)$, 
 with the boundary conditions
\begin{equation}
 D(0)=0,\quad B(N)=0,\quad \mathcal{N}<N,
\end{equation} 
in which the latter two conditions apply only to the finite lattice cases.
The shift operators $e^{\pm\partial}$ act on 
a function defined on the above non-negative integer 
lattice points and 
shift it either $\pm1$:
\begin{equation*}
e^{\pm\partial}\psi(x)=\psi(x\pm 1).
\end{equation*}
The polynomials are assembled into four groups according to 
the form of the sinusoidal coordinate 
$\eta(x)=x,x(x+d)$ and $q^{\pm x}$-linear and $q^{\pm x}$-bilinear. 
In each group, the finite lattice cases are followed by infinite lattice ones. 
The most generic member will be followed by simpler ones. 
For a comprehensive exposition of these polynomials 
in the discrete quantum mechanics formulation and
their applications, see \cite{os12,os24,bdproc}.

%%%%%%%%%%%%%%%%
\subsection{Polynomials having  $\eta(x)=x$, $[0,1,\ldots,N]$ or  
$[0,1,\ldots,\infty)$}
\label{sec:relineta}
This group consists of four polynomials, the Hahn \S\ref{sec:Hahn}, 
the Krawtchouk \S\ref{sec:Kraw} (on  finite lattices)  and the Meixner 
\S\ref{sec:Meix} and the Charlier \S\ref{sec:Char}, 
both on  infinite lattices.
The matrix $\mathcal{M}$ \eqref{Mdef} for this group reads
\begin{align} 
\mathcal{M}_{n\,m}&=\delta_{n\,m}\left(
 \frac{B(x_n)\prod_{j\neq n}^{\mathcal N}(x_n+1-x_j)+
D(x_n)\prod_{j\neq n}^{\mathcal N}(x_n-1-x_j)}
{\prod_{j\neq n}^{\mathcal N}(x_n-x_j)}\right.\n[4pt]
&\hspace{70mm}
\left. +\mathcal{E}(\mathcal{N})-B(x_n)-D(x_n)
\phantom{\frac{(x_n)}{x_n}}\hspace*{-9mm}\right)\n
&\quad +(1-\delta_{n\,m})\frac{\left(B(x_n)\prod_{j\neq n,m}^{\mathcal N}(x_n+1-x_j)-
D(x_n)\prod_{j\neq n,m}^{\mathcal N}(x_n-1-x_j)\right)}
{\prod_{j\neq n}^{\mathcal N}(x_n-x_j)}.
\label{MHahn}
\end{align}
It has eigenvalues $\mathcal{E}(\mathcal{N})-\mathcal{E}(m)$, $m=0,1,
\ldots,\mathcal{N}-1$,
for arbitrary distinct complex values of $\{x_n\}$.
For the zeros $\{x_n\}$ of $P_{\mathcal N}$, $P_{\mathcal N}(x_n)=0$, 
it is  straightforward
to verify Theorem \ref{theo:2} numerically for lower $\mathcal{N}$. 
That is,  the eigenvectors of the matrix 
$\mathcal{M}$ \eqref{MHahn} generate
the lower degree polynomials $\{P_m(x)\}$, $m=0,1,\ldots,\mathcal{N}-1$ 
as in \eqref{Pmgen}.

For these polynomials, the matrix $\mathcal{M}$ \eqref{MHahn} simplifies by using
the algebraic equations for the zeros $\{x_n\}$:
\begin{align}
&B(x_n)\prod_{j\neq n}^{\mathcal N}(x_n+1-x_j)
=D(x_n)\prod_{j\neq n}^{\mathcal N}(x_n-1-x_j),\\
\mathcal{M}_{n\,m}&=\delta_{n\,m}\left(B(x_n)\prod_{j\neq n}\bigl(1+\frac1{x_n-x_j})
+\mathcal{E}(\mathcal{N})-B(x_n)-D(x_n)\right)\n
&+(1-\delta_{n\,m})\,B(x_n)\prod_{j\neq n}
\bigl(1+\frac1{x_n-x_j})\cdot\left(\frac1{x_n+1-x_m}-
\frac1{x_n-1-x_m}\right).
\end{align}

%%%%%%%%%%%%%%%%%%%%%%%%%%%%%%%%%%%%%%%%%%%
%  Hahn                           %
%%%%%%%%%%%%%%%%%%%%%%%%%%%%%%%%%%%%%%%%%%%
\subsubsection{Hahn}
\label{sec:Hahn}
This polynomial depends on two real parameters $a$ 
and $b$ besides $N$ \cite{koekswart,os12}:
\begin{align}
  &B(x)=(x+a)(N-x),\quad
  D(x)= x(b+N-x),\quad a>0, \ b>0,
  \label{hahnBD}\\
 &\phi_0(x)^2
  =\frac{N!}{x!\,(N-x)!}\,\frac{(a)_x\,(b)_{N-x}}{(b)_N}\,,
  \quad  \mathcal{E}(n)= n(n+a+b-1),
  \label{hahnphi0}\\
  &P_n(x)
  ={}_3F_2\Bigl(
  \genfrac{}{}{0pt}{}{-n,\,n+a+b-1,\,-x}
  {a,\,-N}\Bigm|1\Bigr).
  \label{hahnpolydef}
\end{align}

%%%%%%%%%%%%%%%%%%%%%%%%%%%%%%%%%%%%%%%%%%%
%  Krawtchouk                   %
%%%%%%%%%%%%%%%%%%%%%%%%%%%%%%%%%%%%%%%%%%%
\subsubsection{Krawtchouk}
\label{sec:Kraw}
This polynomial depends on one real parameter $p$  
besides $N$ \cite{koekswart,os12}:
\begin{align}
  &B(x)=p(N-x),\quad
  D(x)=(1-p)x,\quad 0<p<1,\quad \mathcal{E}(n)=n,\\
 &\phi_0(x)^2=
  \frac{N!}{x!\,(N-x)!}\Bigl(\frac{p}{1-p}\Bigr)^x,\quad
  P_n(x)
  ={}_2F_1\Bigl(
  \genfrac{}{}{0pt}{}{-n,\,-x}{-N}\Bigm|p^{-1}\Bigr).
  \label{krawtpol}
\end{align}
This polynomial \eqref{krawtpol} is symmetric under the 
interchange $x\leftrightarrow n$, with
$\mathcal{E}(n)=n$ and $\eta(x)=x$.
%%%%%%%%%%%%%%%%%%%%%%%%%%%%%%%%%%%%%%%%%%%
%  Meixner                       %
%%%%%%%%%%%%%%%%%%%%%%%%%%%%%%%%%%%%%%%%%%%
\subsubsection{Meixner}
\label{sec:Meix}

This is another example of the self-dual polynomial \eqref{MeixnerP},
which is symmetric under the interchange $x\leftrightarrow n$, with
$\mathcal{E}(n)=n$ and $\eta(x)=x$.
This polynomial depends on two real parameters $\beta$ 
and $c$, \cite{koekswart,os12}:
\begin{align}
  &B(x)=\frac{c}{1-c}(x+\beta),\quad
  D(x)=\frac{1}{1-c}\,x,\qquad \quad \beta>0,\quad 0<c<1,
  \label{MeixnerBD}\\
  &\phi_0(x)^2=\frac{(\beta)_x\,c^x}{x!}\,,
 \quad P_n(x)
  ={}_2F_1\Bigl(
  \genfrac{}{}{0pt}{}{-n,\,-x}{\beta}\Bigm|1-c^{-1}\Bigr), 
  \quad \mathcal{E}(n)=n.
  \label{MeixnerP}
\end{align}

%%%%%%%%%%%%%%%%%%%%%%%%%%%%%%%%%%%%%%%%%%%
%  Charlier                     %
%%%%%%%%%%%%%%%%%%%%%%%%%%%%%%%%%%%%%%%%%%%
\subsubsection{Charlier}
\label{sec:Char}
The Charlier polynomial depends on one real parameter is $a$ 
and is also self-dual 
$x\leftrightarrow n$ \cite{koekswart,os12}:
\begin{align}
  &B(x)=a,\quad
  D(x)=x, \quad a>0,\quad \mathcal{E}(n)=n,
  \label{charlBD}\\
  &\phi_0(x)^2=\frac{a^x}{x!}\,,\quad P_n(x)
  ={}_2F_0\Bigl(
  \genfrac{}{}{0pt}{}{-n,\,-x}{-}\Bigm|-a^{-1}\Bigr).
    \label{charlP}
\end{align}

%%%%%%%%%%%%%%%%
\subsection{Polynomials having  $\eta(x)$ quadratic in $x$}
\label{sec:rex2}
This group consists of two polynomials, the Racah \S\ref{sec:Rac} 
and the dual Hahn \S\ref{sec:dualH},
both on  finite lattices.
The matrix $\mathcal{M}$ \eqref{Mdef} for this group reads
\begin{align} 
\mathcal{M}_{n\,m}&=\delta_{n\,m}\left(
 \frac{B(x_n)\prod_{j\neq n}^{\mathcal N}\bigl(\eta(x_n+1)-y_j\bigr)+
D(x_n)\prod_{j\neq n}^{\mathcal N}\bigl(\eta(x_n-1)-y_j\bigr)}
{\prod_{j\neq n}^{\mathcal N}(y_n-y_j)}\right.\n[4pt]
&\hspace{70mm}
\left. +\mathcal{E}(\mathcal{N})-B(x_n)-D(x_n)
\phantom{\frac{(x_n)}{x_n}}\hspace*{-9mm}\right)\n
&\quad +(1-\delta_{n\,m})\frac{\dot{\eta}(x_m)}{\dot{\eta}(x_n)}
\frac1{\prod_{j\neq n}^{\mathcal N}(y_n-y_j)}\left(
B(x_n)\!\prod_{j\neq m}^{\mathcal N}\bigl(\eta(x_n+1)-y_j\bigr)\right.\n
&
\left. \hspace{65mm}
+D(x_n)\!\prod_{j\neq m}^{\mathcal N}\bigl(\eta(x_n-1)-y_j\bigr)\right).
\label{MRac}
\end{align}
It has eigenvalues $\mathcal{E}(\mathcal{N})-\mathcal{E}(m)$, 
$m=0,1,\ldots,\mathcal{N}-1$,
for arbitrary distinct complex values of $\{x_n\}$ 
except for the poles of $B(x)$ and $D(x)$.
For the zeros $\{y_n=\eta(x_n)\}$ of $P_{\mathcal N}$, $P_{\mathcal N}(y_n)=0$, 
it is  straightforward
to verify Theorem \ref{theo:2} numerically for lower $\mathcal{N}$.
 That is  the eigenvectors of the matrix 
$\mathcal{M}$ \eqref{MHahn} generate
the lower degree polynomials $\{P_m(\eta)\}$, 
$m=0,1,\ldots,\mathcal{N}-1$ as in \eqref{Pmgen}.

%%%%%%%%%%%%%%%%%%%%%%%%%%%%%%%%%%%%%%%%%%%
%  Racah                           %
%%%%%%%%%%%%%%%%%%%%%%%%%%%%%%%%%%%%%%%%%%%

\subsubsection{Racah}
\label{sec:Rac}
We adopt a  parametrisation \cite{os12} designed to reveal the  
symmetry of the difference Racah equation.
This polynomial has four real parameters, $a$, $b$, $c$ and $d$, 
one of which must be equal to $-N$.
Here we choose the parameter ranges
\begin{equation}
  c=-N,\ d>0,\ a>N+d,\ 0<b<1+d.
\end{equation}
The various data are \cite{koekswart,os12}:
\begin{align}
B(x)
  &=-\frac{(x+a)(x+b)(x+c)(x+d)}{(2x+d)(2x+1+d)},\n
  D(x)
  &=-\frac{(x+d-a)(x+d-b)(x+d-c)x}{(2x-1+d)(2x+d)},\\
  \mathcal{E}(n)&= n(n+\tilde{d}),\qquad \tilde{d}\eqdef a+b+c-d-1,\\
  \eta(x)&=x(x+d),\qquad \dot{\eta}(x)=2x+d,\\
    \phi_0(x)^2&=\frac{(a,b,c,d)_x}{(1+d-a,1+d-b,1+d-c,1)_x}\,
  \frac{2x+d}{d},\\
P_n(\eta(x))
  &={}_4F_3\Bigl(
  \genfrac{}{}{0pt}{}{-n,\,n+\tilde{d},\,-x,\,x+d}
  {a,\,b,\,c}\Bigm|1\Bigr).
\end{align}
All the non-$q$ polynomials in this section can be obtained 
from the Racah by reductions.
%%%%%%%%%%%%%%%%%%%%%%%%%%%%%%%%%%%%%%%%%%%
%  dual Hahn                     %
%%%%%%%%%%%%%%%%%%%%%%%%%%%%%%%%%%%%%%%%%%%
\subsubsection{dual Hahn}
\label{sec:dualH}
We adopt the parametrisation of the dual Hahn polynomial 
so that the duality ($x\leftrightarrow n$) with the Hahn
polynomial \S\ref{sec:Hahn} is obvious. 
Thus the parameters $(a,b)$ are  different from the standard
ones \cite{koekswart} for the dual Hahn
polynomial:
\begin{align}
   &B(x)=\frac{(x+a)(x+a+b-1)(N-x)}
  {(2x-1+a+b)(2x+a+b)},\quad a>0,\ \ b>0,
  \label{dualhahnBD1}\\
  &D(x)=\frac{x(x+b-1)(x+a+b+N-1)}
  {(2x-2+a+b)(2x-1+a+b)},
  \label{dualhahnBD2}\\
  &\mathcal{E}(n)=n,\quad
  \eta(x)= x(x+a+b-1),\quad \dot{\eta}(x)=2x+a+b-1,
  \label{dualhahneeta}\\
 &\phi_0(x)^2
  =\frac{N!}{x!\,(N-x)!}\,
  \frac{(a)_x\,(2x+a+b-1)(a+b)_N}{(b)_x\,(x+a+b-1)_{N+1}\,}\,,
  \label{dualhahnphi0}\\
  &P_n(\eta(x))
  ={}_3F_2\Bigl(
  \genfrac{}{}{0pt}{}{-n,\,x+a+b-1,\,-x}
  {a,\,-N}\Bigm|1\Bigr).
  \label{dualhanpolydef}
 \end{align}

%%%%%%%%%%%%%%%%
\subsection{Polynomials having  $\eta(x)$ linear in $q^{-x}$, 
$[0,1,\ldots,N]$ or  $[0,1,\ldots,\infty)$}
\label{sec:rqlin}
This group of polynomials have the sinusoidal coordinate $\eta(x)=q^{-x}-1$.
Seven polynomials belong to this group; $q$-Hahn \S\ref{sec:qHahn},
 quantum $q$-Krawtchouk \S\ref{sec:qqKraw},
$q$-Krawtchouk \S\ref{sec:qKraw},
affine $q$-Krawtchouk \S\ref{sec:aqKraw}, $q$-Meixner \S\ref{sec:qMeix}, 
Al-Salam-Carlitz II \S\ref{sec:ASC2} and $q$-Charlier \S\ref{sec:qChar}.
As for the naming of the polynomials we follow \cite{koekswart}.
The first four examples are on finite lattices. 
The remaining three are on infinite lattices.
The matrix $\mathcal{M}$ \eqref{Mdef} for this group has essentially 
the same structure and properties as those for the Racah \eqref{MRac}.

%%%%%%%%%%%%%%%%%%%%%%%%%%%%%%%%%%%%%%%%%%%
%  $q$-Hahn                     %
%%%%%%%%%%%%%%%%%%%%%%%%%%%%%%%%%%%%%%%%%%%
\subsubsection{$q$-Hahn}
\label{sec:qHahn}
This polynomial depends on two real parameters $a$ and $b$. 
We choose the parameter range
$0<a<1$ and $0<b<1$. Other data are \cite{os12}:
\begin{align}
  &B(x)=(1-aq^x)(q^{x-N}-1),\qquad
  D(x)=aq^{-1}(1-q^x)(q^{x-N}-b),\\
 &\phi_0(x)^2
  =\frac{(q\,;q)_N}{(q\,;q)_x\,(q\,;q)_{N-x}}\,
  \frac{(a;q)_x\,(b\,;q)_{N-x}}{(b\,;q)_N\,a^x}\,, 
  \quad \mathcal{E}(n)
  =(q^{-n}-1)(1-abq^{n-1}),\\
  &P_n(\eta(x))
  ={}_3\phi_2\Bigl(
  \genfrac{}{}{0pt}{}{q^{-n},\,abq^{n-1},\,q^{-x}}
  {a,\,q^{-N}}\Bigm|q\,;q\Bigr).
\end{align}
%%%%%%%%%%%%%%%%%%%%%%%%%%%%%%%%%%%%%%%%%%%
%  quantum $q$-Krawtchouk       %
%%%%%%%%%%%%%%%%%%%%%%%%%%%%%%%%%%%%%%%%%%%
\subsubsection{quantum $q$-Krawtchouk}
\label{sec:qqKraw}
This polynomial depends on one parameter $p>q^{-N}$:
\begin{align}
  &B(x)=p^{-1}q^x(q^{x-N}-1),\quad
  D(x)=(1-q^x)(1-p^{-1}q^{x-N-1}),\quad \mathcal{E}(n)=1-q^n,\\
  &\phi_0(x)^2
  =\frac{(q\,;q)_N}{(q\,;q)_x(q\,;q)_{N-x}}\,
  \frac{p^{-x}q^{x(x-1-N)}}{(p^{-1}q^{-N}\,;q)_x}\,,\quad
 P_n(\eta(x))
  ={}_2\phi_1\Bigl(
  \genfrac{}{}{0pt}{}{q^{-n},\,q^{-x}}{q^{-N}}\Bigm|q\,;pq^{n+1}\Bigr). 
\end{align}

%%%%%%%%%%%%%%%%%%%%%%%%%%%%%%%%%%%%%%%%%%%%%%%%%%%%%%%%%%%%%
%  $q$-Krawtchouk %
%%%%%%%%%%%%%%%%%%%%%%%%%%%%%%%%%%%%%%%%%%%%%%%%%%%%%%%%%%%%%
\subsubsection{$q$-Krawtchouk}
\label{sec:qKraw}
This polynomial depends on one positive parameter $p>0$ \cite{koekswart,os12}:
\begin{align}
  &B(x)=q^{x-N}-1,\quad
  D(x)=p(1-q^x),\qquad \mathcal{E}(n)=(q^{-n}-1)(1+pq^n), \\
  &\phi_0(x)^2=\frac{(q\,;q)_N}{(q\,;q)_x(q\,;q)_{N-x}}\,
  p^{-x}q^{\frac12x(x-1)-xN},\quad
P_n(\eta(x))
  ={}_3\phi_2\Bigl(
  \genfrac{}{}{0pt}{}{q^{-n},\,q^{-x},\,-pq^n}{q^{-N},\,0}\Bigm|q\,;q\Bigr).
\end{align}

%%%%%%%%%%%%%%%%%%%%%%%%%%%%%%%%%%%%%%%%%%%
%  affine $q$-Krawtchouk         %
%%%%%%%%%%%%%%%%%%%%%%%%%%%%%%%%%%%%%%%%%%%
\subsubsection{affine $q$-Krawtchouk}
\label{sec:aqKraw}
This polynomial has one positive parameter $0<p<q^{-1}$ 
and it is self-dual ($x\leftrightarrow n$):
\begin{align}
  &B(x)=(q^{x-N}-1)(1-pq^{x+1}),\quad
  D(x)=pq^{x-N}(1-q^x),\quad \mathcal{E}(n)=q^{-n}-1,\\
  &\phi_0(x)^2=\frac{(q\,;q)_N}{(q\,;q)_x(q\,;q)_{N-x}}\,
  \frac{(pq\,;q)_x}{(pq)^x}\,,\quad
P_n(\eta(x))
  ={}_3\phi_2\Bigl(
  \genfrac{}{}{0pt}{}{q^{-n},\,q^{-x},\,0}{pq,\,q^{-N}}\Bigm|q\,;q\Bigr).
\end{align}
%%%%%%%%%%%%%%%%%%%%%%%%%%%%%%%%%%%%%%%%%%%
%  $q$-Meixner                  %
%%%%%%%%%%%%%%%%%%%%%%%%%%%%%%%%%%%%%%%%%%%
\subsubsection{$q$-Meixner}
\label{sec:qMeix}
This polynomial has two real parameters $b$ and $c$ with $0<b<q^{-1}$, 
$c>0$ and other data are \cite{koekswart,os12}:
\begin{align}
 &B(x)=cq^x(1-bq^{x+1}),\quad
  D(x)=(1-q^x)(1+bcq^x),\qquad \mathcal{E}(n)=1-q^n,\\
 &\phi_0(x)^2=
  \frac{(bq\,;q)_x}{(q,-bcq\,;q)_x}\,c^xq^{\frac12x(x-1)},\quad
P_n(\eta(x))
  ={}_2\phi_1\Bigl(
  \genfrac{}{}{0pt}{}{q^{-n},\,q^{-x}}{bq}\Bigm|q\,;-c^{-1}q^{n+1}\Bigr). %
\end{align}

%%%%%%%%%%%%%%%%%%%%%%%%%%%%%%%%%%%%%%%%%%%%%%%%%%%%%%%%%%%%
%  Al-Salam-Carlitz II [ %
%%%%%%%%%%%%%%%%%%%%%%%%%%%%%%%%%%%%%%%%%%%%%%%%%%%%%%%%%%%%
\subsubsection{Al-Salam-Carlitz II }
\label{sec:ASC2}

The polynomial depends on one real parameter $0<a<q^{-1}$ \cite{koekswart,os12}:
\begin{align}%
  &B(x)=aq^{2x+1},\quad
  D(x)=(1-q^x)(1-aq^x),\quad \mathcal{E}(n)=1-q^n,\\
    &\phi_0(x)^2=\frac{a^xq^{x^2}}{(q,aq\,;q)_x}\,,\quad
    P_n(\eta(x))
  ={}_2\phi_0\Bigl(
  \genfrac{}{}{0pt}{}{q^{-n},\,q^{-x}}{-}\Bigm|q\,;a^{-1}q^n\Bigr).
  \label{alsalamIInorm}
\end{align}

%%%%%%%%%%%%%%%%%%%%%%%%%%%%%%%%%%%%%%%%%%%
%  $q$-Charlier                %
%%%%%%%%%%%%%%%%%%%%%%%%%%%%%%%%%%%%%%%%%%%
\subsubsection{$q$-Charlier }
\label{sec:qChar}
The polynomial depends on one positive parameter $a>0$ \cite{koekswart,os12}:
\begin{align}
  &B(x)=aq^x,\quad
  D(x)=1-q^x,\qquad
\mathcal{E}(n)=1-q^n,\\
 &\phi_0(x)^2=\frac{a^xq^{\frac12x(x-1)}}{(q\,;q)_x}\,,\quad
P_n(\eta(x))
  ={}_2\phi_1\Bigl(
  \genfrac{}{}{0pt}{}{q^{-n},\,q^{-x}}{0}\Bigm|q\,;-a^{-1}q^{n+1}\Bigr).
  \end{align}

%%%%%%%%%%%%%%%%
\subsection{Polynomials having  $\eta(x)$ linear in $q^{x}$,  $[0,1,\ldots,\infty)$}
\label{sec:rqlin2}
This group of polynomials have the sinusoidal coordinate $\eta(x)=1-q^{x}$.
Three polynomials defined on infinite lattices, little $q$-Jacobi \S\ref{sec:lqJac}, little $q$-Laguerre
\S\ref{sec:lqLag} and 
alternative $q$-Charlier \S\ref{sec:aqChar}, belong to this group.
The matrix $\mathcal{M}$ \eqref{Mdef} for this group has essentially 
the same structure and properties as that for the Racah \eqref{MRac}.

%%%%%%%%%%%%%%%%%%%%%%%%%%%%%%%%%%%%%%%%%%%
%  little $q$-Jacobi             %
%%%%%%%%%%%%%%%%%%%%%%%%%%%%%%%%%%%%%%%%%%%
\subsubsection{little $q$-Jacobi }
\label{sec:lqJac}
This polynomial depends on two positive parameters $0<a,b<q^{-1}$ \cite{koekswart,os12}:
\begin{align}
  &B(x)=a(q^{-x}-bq),\quad
  D(x)=q^{-x}-1,\qquad
\mathcal{E}(n)=(q^{-n}-1)(1-abq^{n+1}),\\
   &\phi_0(x)^2=\frac{(bq\,;q)_x}{(q\,;q)_x}(aq)^x,\quad
  P_n(\eta(x))
  =(-a)^{-n}q^{-\frac12n(n+1)}\frac{(aq\,;q)_n}{(bq\,;q)_n}\,
  {}_2\phi_1\Bigl(
  \genfrac{}{}{0pt}{}{q^{-n},\,abq^{n+1}}{aq}\Bigm|q\,;q^{x+1}\Bigr).
\end{align}

%%%%%%%%%%%%%%%%%%%%%%%%%%%%%%%%%%%%%%%%%%%%%%%%%%%%%%%%%%%%%%%%%
%  little $q$-Laguerre/Wall  %
%%%%%%%%%%%%%%%%%%%%%%%%%%%%%%%%%%%%%%%%%%%%%%%%%%%%%%%%%%%%%%%%%
\subsubsection{little $q$-Laguerre/Wall }
\label{sec:lqLag}
Putting $b=0$ in the little $q$-Jacobi \S\ref{sec:lqJac} 
gives this polynomial \cite{koekswart,os12}:
\begin{align}
  &B(x)=aq^{-x},\quad
  D(x)=q^{-x}-1,\qquad
\mathcal{E}(n)=q^{-n}-1,\\
&\phi_0(x)^2=\frac{(aq)^x}{(q\,;q)_x}\,,\quad
  P_n(\eta(x))
  ={}_2\phi_0\Bigl(
  \genfrac{}{}{0pt}{}{q^{-n},\,q^{-x}}{-}\Bigm|q\,;a^{-1}q^x\Bigr).
    \end{align}
%%%%%%%%%%%%%%%%%%%%%%%%%%%%%%%%%%%%%%%%%%%
%  alternative $q$-Charlier    %
%%%%%%%%%%%%%%%%%%%%%%%%%%%%%%%%%%%%%%%%%%%
\subsubsection{alternative $q$-Charlier}
\label{sec:aqChar}
This polynomial depends on one positive parameter $a>0$ \cite{koekswart,os12}:
\begin{align}
  &B(x)=a,\quad
  D(x)=q^{-x}-1,\quad
\mathcal{E}(n)=(q^{-n}-1)(1+aq^n),\\
  &\phi_0(x)^2=\frac{a^xq^{\frac12x(x+1)}}{(q\,;q)_x}\,,\quad
P_n(\eta(x))
  =q^{nx}\,{}_2\phi_1\Bigl(
  \genfrac{}{}{0pt}{}{q^{-n},\,q^{-x}}{0}\Bigm|q\,;-a^{-1}q^{-n+1}\Bigr).
  \end{align}

%%%%%%%%%%%%%%%%
\subsection{Polynomials having  $\eta(x)$ bilinear in $q^{-x}$ and $q^{x}$,  $[0,1,\ldots,N]$}
\label{sec:rqblin}
This group of polynomials have the sinusoidal coordinate $\eta(x)=(q^{-x}-1)(1-Aq^x)$ 
with some constant $A$ and on finite lattices. 
The $q$-Racah \S\ref{sec:qRac} and dual $q$-Hahn \S\ref{sec:dqHahn} belong to this group.
The matrix $\mathcal{M}$ \eqref{Mdef} for this group has essentially the same structure and properties as those for the Racah \eqref{MRac}.

%%%%%%%%%%%%%%%%%%%%%%%%%%%%%%%%%%%%%%%%%%%
%  $q$-Racah                    %
%%%%%%%%%%%%%%%%%%%%%%%%%%%%%%%%%%%%%%%%%%%
\subsubsection{$q$-Racah}
\label{sec:qRac}

 This polynomial is the most
general of all the classical orthogonal polynomials of a discrete variable. 
The rest of the polynomials in \S\ref{sec:dQMre} is obtained 
from $q$-Racah by reductions.
We adopt a  parametrisation \cite{os12} designed to reveal the symmetry of the difference 
$q$-Racah equation.
This polynomial has four real parameters, $a$, $b$, $c$ and $d$, 
one of which must be equal to $q^{-N}$.
The set of parameters  is different from
the standard one $(\alpha,\beta,\gamma,\delta)$ in the same manner
as for the Racah polynomial \S\ref{sec:Rac}.

Here we choose the parameter ranges and introduce $\tilde{d}$:
\begin{equation}
  c=q^{-N},\ 0<d<1,\ 0<a<q^Nd,\ qd<b<1,\quad \tilde{d}\eqdef abcd^{-1}q^{-1}.
\end{equation}
Other data are  \cite{koekswart,os12}:
\begin{align}
  &B(x)
  =-\frac{(1-aq^x)(1-bq^x)(1-cq^x)(1-dq^x)}
  {(1-dq^{2x})(1-dq^{2x+1})}\,,\\
  &D(x)
  =- \tilde{d}\,
  \frac{(1-a^{-1}dq^x)(1-b^{-1}dq^x)(1-c^{-1}dq^x)(1-q^x)}
  {(1-dq^{2x-1})(1-dq^{2x})},\\
  &\mathcal{E}(n)=(q^{-n}-1)(1-\tilde{d}q^n),\quad
  \eta(x)=(q^{-x}-1)(1-dq^x),\\
  &\phi_0(x)^2=\frac{(a,b,c,d\,;q)_x}
  {(a^{-1}dq,b^{-1}dq,c^{-1}dq,q\,;q)_x\,\tilde{d}^x}\,
  \frac{1-dq^{2x}}{1-d}\,,\\
  &P_n(\eta(x))
  ={}_4\phi_3\Bigl(
  \genfrac{}{}{0pt}{}{q^{-n},\,\tilde{d}q^n,\,q^{-x},\,dq^x}
  {a,\,b,\,c}\Bigm|q\,;q\Bigr).
\end{align}

%%%%%%%%%%%%%%%%%%%%%%%%%%%%%%%%%%%%%%%%%%%
%  dual $q$-Hahn                 %
%%%%%%%%%%%%%%%%%%%%%%%%%%%%%%%%%%%%%%%%%%%
\subsubsection{dual $q$-Hahn}
\label{sec:dqHahn}

We adopt the same parameters $(a,b)$ ($0<a,b<1$) for
the $q$-Hahn \S\ref{sec:qHahn} and dual $q$-Hahn polynomials \cite{koekswart,os12}:
\begin{align}
  &B(x)=
  \frac{(q^{x-N}-1)(1-aq^x)(1-abq^{x-1})}
  {(1-abq^{2x-1})(1-abq^{2x})},\\
  &D(x)=aq^{x-N-1}
  \frac{(1-q^x)(1-abq^{x+N-1})(1-bq^{x-1})}
  {(1-abq^{2x-2})(1-abq^{2x-1})},\\
  &\mathcal{E}(n)=q^{-n}-1,\quad
  \eta(x)=(q^{-x}-1)(1-abq^{x-1}),\\
    &\phi_0(x)^2
  =\frac{(q\,;q)_N}{(q\,;q)_x\,(q\,;q)_{N-x}}\,
  \frac{(a,abq^{-1}\,;q)_x}{(abq^N,b\,;q)_x\,a^x}\,
  \frac{1-abq^{2x-1}}{1-abq^{-1}}\,,\\
  &P_n(\eta(x))
  ={}_3\phi_2\Bigl(
  \genfrac{}{}{0pt}{}{q^{-n},\,abq^{x-1},\,q^{-x}}
  {a,\,q^{-N}}\Bigm|q\,;q\Bigr).
\end{align}

%%%%%%%%%%%%%%%%%%%%%%%%%%%%%%%%%%%%%%%%%%%%%%%%%%%%%%%%%%%%%%%
%                                                             %
%  7. Summary and Comments                                    %
%                                                             %
%%%%%%%%%%%%%%%%%%%%%%%%%%%%%%%%%%%%%%%%%%%%%%%%%%%%%%%%%%%%%%%
\section{Summary and Comments}
\label{sec:sum}

 Based on the quantum mechanical reformulation of classical orthogonal polynomials
 \cite{os12,os13,os24}, an $\mathcal{N}\times \mathcal{N}$ matrix $\mathcal{M}$ \eqref{Mdef} describing the small oscillations around the zeros of the degree $\mathcal{N}$
 polynomial is derived.
By construction, its components depend explicitly on the zeros, but its eigenvalues 
\eqref{eigvals} are independent of them. Its eigenvalues are the difference of those of the
differential/difference operator $\widetilde{\mathcal H}$ \eqref{Htil},
which governs the classical orthogonal polynomial, corresponding to the degree $\mathcal{N}$ and a lower degree, Theorem\,\ref{theo:1}.
The corresponding eigenvectors \eqref{eigveceq} of $\mathcal{M}$ provide the
representations of the lower degree polynomials in terms of the zeros of the degree
$\mathcal{N}$ polynomial, \eqref{Pmgen}, Theorem\,\ref{theo:2}.
It should be stressed that these theorems are valid {\em universally for all the classical orthogonal polynomials\/}.
We have provided the necessary data for most of the classical orthogonal polynomials
ranging from the Hermite, Laguerre, Jacobi, Wilson, Askey-Wilson, Racah, $q$-Racah
and their reduced form polynomials for self contained verification of the main results.
The data include the proper ranges of the parameters, which are important for numerical verification.

The ingredients of the matrix $\mathcal{M}$ \eqref{Mdef} are the {\em sinusoidal coordinate\/} $\eta(x)$ \eqref{fac} and the differential/difference
operator $\widetilde{\mathcal H}$, that is, the analytic functions $V(x)$ and $V^*(x)$ 
\eqref{HVV*} for the Wilson, Askey-Wilson and their reduced form polynomials,
and the two non-negative functions $B(x)$ and $D(x)$ \eqref{realHt} for the Racah, $q$-Racah and their reduced form polynomials.
The close relationship between these functions and the sinusoidal coordinates was
elucidated in \cite{os14}.

The present research was inspired by a recent work of Bihun and Calogero \cite{bihun-cal}
discussing the properties of the zeros of the polynomials belonging to the Askey scheme.
Their paper is a certain generalisation of the old results by Calogero and his co-authors 
\cite{ahmed} on the properties of the zeros of the Hermite, Laguerre and Jacobi polynomials.
Although our matrix $\mathcal{M}$ \eqref{Mdef} and their matrices in \cite{bihun-cal,ahmed}
have related eigenvalues, our $\mathcal{M}$ \eqref{Mdef}  is conceptually and structurally different from those matrices in the two papers \cite{bihun-cal, ahmed}.
The motivation and guiding principle of the earlier works 
\cite{stiel,szego,bihun-cal, ahmed}, \cite{calmat}--\cite{os6} were the Diophantine properties of the Hessian matrices
describing the small oscillations around the equilibrium of {\em exactly solvable
multi-particle dynamics\/}.
In contrast, our matrix $\mathcal{M}$ \eqref{Mdef} describes the perturbations around the
zeros of a classical orthogonal polynomial of a {\em single variable\/}.

For the Classical orthogonal polynomials, {\em i.e.\/} the Hermite, Laguerre and Jacobi
polynomials, we have demonstrated algebraic derivation of the Theorems in \S\ref{sec:oQM} with the help of the old results in \cite{ahmed}.
These contents can be reformulated by using `finite dimensional representations of differential operators', which was developed by Calogero \cite{book}.
It is a good challenge to deliver similar algebraic derivation of the {\bf Corollary}\,\ref{coro:1} \eqref{coro1eq} for each of the polynomials 
in the Askey scheme.

After completing this paper, we noticed a recent publication \cite{bihun-cal2}
discussing related matrices for $q$-Askey scheme polynomials and another discussing finite dimensional representations of difference operators 
\cite{fincal}. The latter might be useful for algebraic derivation of the {\bf Corollary}\,\ref{coro:1} \eqref{coro1eq} for the ($q$-)Askey scheme
of hypergeometric orthogonal polynomials. 
%%%%%%%%%%%%%%%%%%%%%%%%%%%%%%%%%%%%%%%%%%%%%%%%%%%%%%%%%%%%%%%
%                                                             %
%  Acknowledgments                                            %
%                                                             %
%%%%%%%%%%%%%%%%%%%%%%%%%%%%%%%%%%%%%%%%%%%%%%%%%%%%%%%%%%%%%%%
\section*{Acknowledgements}
It is a pleasure to thank the organizers of the CRM-ICMAT Workshop on 
``Exceptional orthogonal polynomials and exact solutions in mathematical physics" 
(Segovia, Spain, 7--12 July 2014). 
R.\,S. thanks Francesco Calogero and Kazuhiko Aomoto for useful and insightful discussion and comments.
He thanks Pauchy Hwang and Department of Physics, National Taiwan University
for hospitality.

%%%%%%%%%%%%%%%%%%%%%%%%%%%%%%%%%%%%%%%%%%%%%%%%%%%%%%%%%%%%%%%
%                                                             %
%  Appendix                                               %
%                                                             %
%%%%%%%%%%%%%%%%%%%%%%%%%%%%%%%%%%%%%%%%%%%%%%%%%%%%%%%%%%%%%%%
\section*{Appendix: Symbols and definitions}
\label{append}
\setcounter{equation}{0}
\renewcommand{\theequation}{A.\arabic{equation}}

For self-containedness we collect several definitions related to
the ($q$-)hypergeometric functions \cite{koekswart}.

\noindent
%%%%%%%%%%%%%%%%%%%
$\circ$ Shifted factorial $(a)_n$ :
\begin{equation}
  (a)_n\eqdef\prod_{k=1}^n(a+k-1)=a(a+1)\cdots(a+n-1)
  =\frac{\Gamma(a+n)}{\Gamma(a)}.
  \label{defPoch}
\end{equation}
%%%%%%%%%%%%%%%%%%%
$\circ$ $q$-Shifted factorial $(a\,;q)_n$ :
\begin{equation}
  (a\,;q)_n\eqdef\prod_{k=1}^n(1-aq^{k-1})=(1-a)(1-aq)\cdots(1-aq^{n-1}).
  \label{defqPoch}
\end{equation}
%%%%%%%%%%%%%%%%%%%
$\circ$ hypergeometric functions ${}_rF_s$ :
\begin{equation}
  {}_rF_s\Bigl(\genfrac{}{}{0pt}{}{a_1,\,\cdots,a_r}{b_1,\,\cdots,b_s}
  \Bigm|z\Bigr)
  \eqdef\sum_{n=0}^{\infty}
  \frac{(a_1,\,\cdots,a_r)_n}{(b_1,\,\cdots,b_s)_n}\frac{z^n}{n!}\,,
  \label{defhypergeom}
\end{equation}
where $(a_1,\,\cdots,a_r)_n\eqdef\prod_{j=1}^r(a_j)_n
=(a_1)_n\cdots(a_r)_n$.\\
%%%%%%%%%%%%%%%%%%%
$\circ$ $q$-hypergeometric functions (the basic hypergeometric functions)
${}_r\phi_s$ :
\begin{equation}
  {}_r\phi_s\Bigl(
  \genfrac{}{}{0pt}{}{a_1,\,\cdots,a_r}{b_1,\,\cdots,b_s}
  \Bigm|q\,;z\Bigr)
  \eqdef\sum_{n=0}^{\infty}
  \frac{(a_1,\,\cdots,a_r\,;q)_n}{(b_1,\,\cdots,b_s\,;q)_n}
  (-1)^{(1+s-r)n}q^{(1+s-r)n(n-1)/2}\frac{z^n}{(q\,;q)_n}\,,
  \label{defqhypergeom}
\end{equation}
where $(a_1,\,\cdots,a_r\,;q)_n\eqdef\prod_{j=1}^r(a_j\,;q)_n
=(a_1\,;q)_n\cdots(a_r\,;q)_n$.

%%%%%%%%%%%%%%%%%%%%%%%%%%%%%%%%%%%%%%%%%%%%%%%%%%%%%%%%%%%%%%%
%                                                             %
%  References                                                 %
%                                                             %
%%%%%%%%%%%%%%%%%%%%%%%%%%%%%%%%%%%%%%%%%%%%%%%%%%%%%%%%%%%%%%%

%\goodbreak


\begin{thebibliography}{99}
% 
% for hyphenation : \hspace{0pt}

\bibitem{stiel}
T.\,J.\, Stieltjes, 
{\it Ouvres Compl\'etes\/}, vol.2 Noordhoff, Groningen (1918).

\bibitem{szego}
 G.\,Szeg\"{o},
{\it Orthogonal polynomials\/}, Fourth edition,
Amer. Math. Soc. New York (1975).

\bibitem{chihara}
T.\,S.\,Chihara,
{\it An Introduction to orthogonal polynomials},
Gordon and Breach, New York (1978).

\bibitem{askey}
G.\,E.\,Andrews, R.\,Askey and R.\,Roy,
{\it Special Functions},
Encyclopedia of mathematics and its applications,
Cambridge Univ. Press, Cambridge (1999).


\bibitem{ismail}
M.\,E.\,H.\,Ismail,
{\it Classical and quantum orthogonal polynomials in one variable},
Encyclopedia of mathematics and its applications,
Cambridge Univ. Press, Cambridge (2005).

\bibitem{koekswart}
R.\,Koekoek and R.\,F.\,Swarttouw,
``The Askey-scheme of hypergeometric orthogonal polynomials and
its $q$-analogue,''
{\tt arXiv:math.CA/9602214};
R.\,Koekoek, P.\,A.\,Lesky and R.\,F.\,\hspace{0pt}Swarttouw, 
{\it Hypergeometric orthogonal polynomials and their $q$-analogues,\/}
Springer-Verlag (2010).

%%%%%%% multi-index
\bibitem{os25}
S.\,Odake and R.\,Sasaki,
``Exactly solvable quantum mechanics and infinite families of
multi-indexed orthogonal polynomials,"
Phys. Lett. {\bf B702} (2011) 164-170,
{\tt arXiv:\hspace{0pt}1105.0508[math-ph]};
``Multi-indexed ($q$-)Racah polynomials,"
J. Phys. {\bf A45} (2012) 385201 (21 pp),
{\tt arXiv:1203.5868[math-ph]}; %, doi:10.1088/1751-8113/45/38/385201;
``Multi-indexed Wilson and Askey-Wilson polynomials,"
J. Phys. {\bf A46} (2013) 045204 (22 pp),
{\tt arXiv:1207.5584[math-ph]}.%, doi:10.1088/1751-8113/46/4/045204. 

%%%%%%%%%xop
\bibitem{xop}
D.\,G\'{o}mez-Ullate, N.\,Kamran and R.\,Milson,
``An extended class of orthogonal polynomials defined by a Sturm-Liouville problem,''
J. Math. Anal. Appl. {\bf 359} (2009) 352-367,
{\tt arXiv:0807.3939\hspace{0pt}[math-ph]};
C.\,Quesne,
``Exceptional orthogonal polynomials, exactly solvable potentials
and supersymmetry,''
J. Phys. {\bf A41} (2008) 392001,
{\tt arXiv:0807.4087[quant-ph]};
S.\,Odake and R.\,Sasaki,
``Infinitely many shape invariant potentials and new orthogonal
polynomials,''
Phys. Lett. {\bf B679} (2009) 414-417,
{\tt arXiv:0906.0142[math-ph]}.

\bibitem{bihun-cal}
O. Bihun and F. Calogero, 
``Properties of the zeros of the polynomials belonging to the Askey scheme,"
{\tt arXiv:1407.3379[math.CA].}

\bibitem{infhull}
L.\,Infeld and T.\,E.\,Hull,
``The factorization method,''
Rev. Mod. Phys. {\bf 23} (1951) 21-68.

\bibitem{susyqm}
F.\,Cooper, A.\,Khare and U.\,Sukhatme,
``Supersymmetry and quantum mechanics,''
Phys. Rep. {\bf 251} (1995) 267-385.



\bibitem{os12}
S.\,Odake and R.\,Sasaki,
``Orthogonal Polynomials from Hermitian Matrices,"
J. Math. Phys. {\bf 49} (2008) 053503 (43pp),
{\tt arXiv:0712.4106[math.CA]}.

\bibitem{os13}
S.\,Odake and R.\,Sasaki,
``Exactly solvable `discrete' quantum mechanics;
shape invariance, Heisenberg solutions,
annihilation-creation operators and coherent states,"
Prog. Theor. Phys. {\bf 119} (2008) 663-700,
{\tt arXiv:0802.1075[quant-ph]}.

\bibitem{os24}
S.\,Odake and R.\,Sasaki,
``Discrete quantum mechanics," (Topical Review)
J. Phys. {\bf A44} (2011) 353001 (47 pp),
{\tt arXiv:1104.0473[math-ph]}.
%%%%%  -------------------- 


%%%%%%% sinusoidal
\bibitem{os7}
S.\,Odake and R.\,Sasaki,
``Unified theory of annihilation-creation operators for solvable
(`discrete') quantum mechanics,''
J. Math. Phys. {\bf 47} (2006) 102102 (33pp),
{\tt arXiv:\hspace{0pt}quant-ph/0605215};
``Exact solution in the Heisenberg picture and annihilation-creation
operators,"
Phys. Lett. {\bf B641} (2006) 112-117,
{\tt arXiv:quant-ph/0605221}.


\bibitem{ahmed}
S.\, Ahmed, M.\,Bruschi, F. Calogero, M.\,A.\, Olshanetsky and A.\, Perelomov,
``Properties of the zeros of the Classical polynomials and of the Bessel functions,"
Nuouvo Cimento {\bf 49} (1979) 173-199.






\bibitem{os11}
S.\,Odake and R.\,Sasaki,
``$q$-oscillator from the $q$-Hermite Polynomial,''
Phys. Lett. {\bf B663} (2008) 141-145,
{\tt arXiv:0710.\hspace{0pt}2209[hep-th]}.


\bibitem{nikiforov}
A.\,F.\,Nikiforov, S.\,K.\,Suslov, and V.\,B.\,Uvarov,
{\it Classical Orthogonal Polynomials of a Discrete Variable\/},
Springer, Berlin (1991).



\bibitem{bdproc}
R.\,Sasaki,
``Exactly solvable birth and death processes,"
J. Math. Phys. {\bf 50} (2009) 103509 (18pp),
{\tt arXiv:0903.3097[math-ph]}.

\bibitem{os14}
S.\,Odake and R.\,Sasaki,
``Unified theory of exactly and quasi-exactly solvable `discrete'
quantum mechanics: I. Formalism,"
J. Math. Phys {\bf 51} (2010) 083502 (24pp),
{\tt arXiv:0903.2604[math-ph]}.

\bibitem{calmat}
F.\,Calogero, ``On the zeros of the classical polynomials,'' Lett. Nuovo
Cim. {\bf 19} (1977) 505-507;
``Equilibrium configuration of one-dimensional many-body problems
with quadratic and inverse quadratic pair potentials,"
Lett. Nuovo Cim. {\bf 22} (1977) 251-253;
``Eigenvectors of a matrix related to the zeros of Hermite polynomials,"
Lett. Nuovo Cim. {\bf 24} (1979) 601-604;
``Matrices, differential operators and polynomials'',
J. Math. Phys. {\bf 22} (1981) 919-934.

\bibitem{CorSa}
E.\,Corrigan and R.\,Sasaki,
``Quantum vs Classical  Integrability in Calogero-Moser
Systems",
J. Phys. {\bf A35}  (2002) 7017-7061, {\tt arXiv:hep-th/0204039}.

\bibitem{RagSa}
O. Ragnisco and R.\,Sasaki,
``Quantum vs Classical Integrability in Ruijsenaars-Schneider Systems,''
J. Phys.   {\bf A37} (2004) 469 - 479,  {\tt arXiv:hep-th/0305120}.  

\bibitem{os2}
S.~Odake and R.~Sasaki,
``Equilibria of 'discrete' integrable systems and deformations of classical
orthogonal polynomials,''
J. Phys. {\bf A37} (2004) 11841-11876, {\tt arXiv:hep-th/0407155}.

\bibitem{vd}
J.\,F.\, van Diejen,
``On the Equilibrium Configuration of the $BC$-type Ruijsenaars-Schneider
System,"
J. Nonlinear Math. Phys. {\bf 12} Suppl. 1 (2005) 689-696,
{\tt arXiv:\hspace{0pt}math-ph\hspace{0pt}/0410008}.

\bibitem{os6}
S. Odake and R. Sasaki,
``Equilibrium Positions, Shape Invariance and Askey-Wilson
       Polynomials,"
J. Math. Phys. {\bf 46} (2005)  063513, 10 pages, {\tt arXiv:hep-th/0410109}.

\bibitem{book}
F. Calogero, {\it  Classical many-body problems amenable to exact treatments}, Lecture Notes in Physics Monograph {\bf m66}, Springer, Berlin (2001).

\bibitem{bihun-cal2}
O. Bihun and F. Calogero, 
``Properties of the zeros of the polynomials
belonging to the q-Askey scheme?," {\tt arXiv:1410.0549[math-ph]}.

\bibitem{fincal}
F. Calogero, 
``Finite-dimensional representations of difference
operators, and the identification of some remarkable matrices," (to be published).
\end{thebibliography}
\end{document}